# ON THE SHINTANI ZETA FUNCTION FOR THE SPACE OF PAIRS OF BINARY HERMITIAN FORMS

Akihiko Yukie[1]

Oklahoma State University

## Introduction

Throughout this paper, $k$ is a number field. We fix a quadratic extension $k_1 = k(\alpha_0)$ of $k$, where $\alpha_0 = \sqrt{\beta_0}$ for a certain $\beta_0 \in k^\times \setminus (k^\times)^2$. In this paper, we consider the zeta function defined for the space of pairs of binary Hermitian forms. This is the prehomogeneous vector space we discussed in §2 [1] and is a non-split form of the $D_4$ case in [4].

The purpose of this paper is to determine the principal part of the adjusted zeta function. Our main result is Theorem (8.15). Our case resembles the space of pairs of binary quadratic forms which we discussed in Chapter 5 [5]. However, the meaning of the problem, the pole structure of the adjusted zeta function, and the adjusting terms are different from those of the above case. The zeta function of our case is a counting function of the class number times the regulator of fields of the form $k(\sqrt{\beta_0}, \sqrt{\beta})$ with $\beta_0$ fixed and $\beta \in k$, whereas the zeta function for the above case is essentially the same as the space of binary quadratic forms. Our case is the first example where we need two adjusting terms whereas the adjusting term of the above case was essentially the same as that of the space of binary quadratic forms. We will not prove the meromorphic continuation of the adjusting terms in this paper. For the sake of the density theorem, we still get the residue of the rightmost pole of the zeta function (without adjusting). If this space appears in a bigger prehomogeneous vector space, we expect to use special values of the adjusted zeta function to describe the principal part of the zeta function of the bigger prehomogeneous vector space. So the principal part formula of the adjusted zeta function should be enough for most applications. We have yet to obtain the corresponding density theorem because we have not carried out the local theory.

In §1, we define the prehomogeneous vector space in this paper, and prove a stratification of this vector space. Since our group is non-split, we cannot directly apply Kempf's theorem on the rationality of the equivariant Morse stratification [2], even though it does not seem so difficult to make necessary adjustments. However, since our representation is rather small, we will handle the stratification explicitly. We will probably consider such adjustments for general non-split prehomogeneous

[1]The author is partially supported by NSF grant DMS-9401391 and NSA grant MDA-904-93-H-3035

vector spaces in the future. In §2, we define the zeta function, the Fourier transform, and the smoothed Eisenstein series. Our method for computing the principal part of the zeta function is similar to that of [5]. In §3, we review the pole structure of the zeta function for the space of (single) binary Hermitian forms. In §§4,5, we consider the adjusting terms. In §§6,7, we consider contributions from unstable points. In §8, we determine the principal part of the adjusted zeta function.

For basic notions on adeles, see [3]. The ring of adeles, the group of ideles, and the adelic absolute value of $k$ are denoted by $\mathbb{A}$, $\mathbb{A}^\times$, $|\ |$ respectively. We fix a non-trivial character $<\ \ >$ of $\mathbb{A}/k$. Note that $< \mathrm{Tr}_{k_1/k}(\ \ ) >$ defines a non-trivial character of $\mathbb{A}_{k_1}/k_1$. The ring of adeles, the group of ideles, and the adelic absolute value of $k_1$ are denoted by $\mathbb{A}_{k_1}$, $\mathbb{A}_{k_1}^\times$, $|\ |_{k_1}$ respectively. Note that by the inclusion $\mathbb{A} \to \mathbb{A}_{k_1}$, an idele $(a_v)_v$ corresponds to the idele $(b_w)_w$ such that $b_w = a_v$ if $w$ is a place over $v$. Let $\mathbb{A}^1 = \{t \in \mathbb{A}^\times \mid |t| = 1\}$, $\mathbb{A}_{k_1}^1 = \{t \in \mathbb{A}_{k_1}^\times \mid |t|_{k_1} = 1\}$. Identifying $k_1 \otimes \mathbb{A} \cong \mathbb{A}_{k_1}$ the norm map $\mathrm{N}_{k_1/k}$ can be extended to a map from $\mathbb{A}_{k_1}$ to $\mathbb{A}$. It is known (see [3, p. 139]) that $|\mathrm{N}_{k_1/k}(t)| = |t|_{k_1}$ for $t \in \mathbb{A}_{k_1}$. The set of positive real numbers is denoted by $\mathbb{R}_+$.

Suppose $[k : \mathbb{Q}] = n$. Then $[k_1 : \mathbb{Q}] = 2n$. For $\lambda \in \mathbb{R}_+$, $\underline{\lambda} \in \mathbb{A}^\times$ is the idele whose component at any infinite place is $\lambda^{\frac{1}{n}}$ and whose component at any finite place is 1. Also $\underline{\lambda}_{k_1} \in \mathbb{A}_{k_1}^\times$ is the idele whose component at any infinite place is $\lambda^{\frac{1}{2n}}$ and whose component at any finite place is 1. Clearly, $\underline{\lambda} = \underline{\lambda}_{k_1}^2$. Since $|\underline{\lambda}| = \lambda$, $|\underline{\lambda}_{k_1}|_{k_1} = \lambda$, this means that $|\underline{\lambda}|_{k_1} = \lambda^2$.

If $X$ is a variety over $k$ and $R$ is a $k$–algebra, the set of $R$–rational points of $X$ is denoted by $X_R$. The space of Schwartz–Bruhat functions on $V_{\mathbb{A}}$ is denoted by $\mathscr{S}(V_{\mathbb{A}})$. If $f, g$ are functions on a set $X$, $f \ll g$ means that there is a constant $C$ such that $f(x) \le Cg(x)$ for all $x \in X$.

## §1 The space of pairs of binary Hermitian forms

Let $G_1 = \mathrm{GL}(2)_{k_1}$, $G_2 = \mathrm{GL}(2)_k$, and $G = G_1 \times G_2$. We consider $G$ as a group over $k$. Let $\mathrm{M}(2,2)_{k_1}$ be the set of $2 \times 2$ matrices whose entries are in $k_1$. Let $\sigma$ be the non-trivial element of $\mathrm{Gal}(k_1/k)$. An element $x \in \mathrm{M}(2,2)_{k_1}$ is called a binary Hermitian form if ${}^t x = x^\sigma$. Let $W$ be the space of binary Hermitian forms and $V = W \oplus W$. We identify $x = (x_1, x_2) \in V$ with $M_x(v) = v_1 x_1 + v_2 x_2$ which is a $2 \times 2$ matrix with entries in linear forms in two variables $v = (v_1, v_2)$. We define an action of $G$ on $V$ by

$$(g_1, g_2) M_x(v) = g_1 M_x(v g_2) {}^t g_1^\sigma.$$

We define $F_x(v) = \det M_x(v)$. We proved in [1] that $(G, V)$ is a prehomogeneous vector space and the relative invariant polynomial is the discriminant of $F_x(v)$.

We can express

$$x_1 = \begin{pmatrix} x_{10} & x_{11} \\ x_{11}^\sigma & x_{12} \end{pmatrix},\ x_2 = \begin{pmatrix} x_{20} & x_{21} \\ x_{21}^\sigma & x_{22} \end{pmatrix},$$

where $x_{10}, x_{12}, x_{20}, x_{22} \in k$, $x_{11}, x_{21} \in k_1$, or $x_{10}, x_{12}, x_{20}, x_{22} \in \mathbb{A}$, $x_{11}, x_{21} \in \mathbb{A}_{k_1}$. We choose $x = (x_{ij})$ $(i = 1, 2,\ j = 0, 1, 2)$ as the coordinate system for $V$. Note that we are considering $k_1$ as a two dimensional vector space over $k$ and we can write $x_{11} = x_{11,1} + x_{11,2}\alpha_0$ where $x_{11,1}, x_{11,2} \in k$. But it turns out that $x_{11,1}$ and $x_{11,2}$ have the same weight with respect to a maximal split torus of $G$. So for our purposes, we don't have to use any $k$–coordinate of $x_{11}$ and we leave $x_{11}$ as an element of $k_1$. The situation is the same for $x_{21}$.

For $i = 1, 2$, let $T_i \subset G_i$ be the subgroup of diagonal matrices, $N_i \subset G_i$ the subgroup of lower triangular matrices with diagonal entries 1, and $B_i = T_i N_i$. We define $T = T_1 \times T_2$, $N = N_1 \times N_2$. Then $B = TN$ is a Borel subgroup of $G$. Let

$$\nu = \begin{pmatrix} 0 & 1 \\ 1 & 0 \end{pmatrix}.$$

The Weyl group of $G$ is generated by $(\nu, 1)$ and $(1, \nu)$.

For the rest of this section, we will consider a stratification of $V_k$. Let

$$\begin{aligned}
&Y_1 = \{x \in V \mid x_{10} = x_{11} = x_{20} = 0\},\ Z_1 = \{x \in Y_1 \mid x_{22} = 0\},\\
&Y_1^{\mathrm{ss}} = \{x \in Y_1 \mid x_{12}, x_{21} \neq 0\},\ Z_1^{\mathrm{ss}} = Y_1^{\mathrm{ss}} \cap Z_1,\\
&Y_2 = Z_2 = \{x \in V \mid x_{1i} = 0 \text{ for } i = 0, 1, 2\},\\
&Y_2^{\mathrm{ss}} = Z_2^{\mathrm{ss}} = \{x \in Y_2 \mid \det x_2 \neq 0\},\\
&Y_3 = Z_3 = \{x \in V \mid x_{1i} = 0 \text{ for } i = 0, 1, 2, x_{20} = x_{21} = 0\},\\
&Y_3^{\mathrm{ss}} = Z_3^{\mathrm{ss}} = \{x \in Z_3 \mid x_{22} \neq 0\}.
\end{aligned}$$

We define $S_i = GY_i^{\mathrm{ss}}$ for $i = 1, 2, 3$. Let $P_2 = G_1 \times B_2$, $P_1 = P_3 = B$.

**Proposition (1.1)** (1) $V_k \setminus \{0\} = V_k^{\mathrm{ss}} \coprod S_{1k} \coprod S_{2k} \coprod S_{3k}$.
(2) $S_{ik} \cong G_k \times_{P_{ik}} Y_{ik}^{\mathrm{ss}}$ for $i = 1, 2$.

*Proof.* We first prove that $V_k \setminus \{0\} = V_k^{\mathrm{ss}} \cup S_{1k} \cup S_{2k} \cup S_{3k}$. By an easy computation,

$$F_x(v) = \det x_1 v_1^2 + Q(x) v_1 v_2 + \det x_2 v_2^2,$$

where $Q(x) = x_{10}x_{22} + x_{12}x_{20} - \mathrm{Tr}_{k_1/k}(x_{11}x_{21}^{\sigma})$. A point $x \in V_k \setminus \{0\}$ is an unstable point if and only if the above polynomial of $v = (v_1, v_2)$ either has a square factor or is zero. If a rational quadratic polynomial has a square factor, that factor is rational. So this means that $x$ is unstable if and only if there exists $g_2 \in G_{2k}$ such that if $y = (y_1, y_2) = g_2 x$, $\det y_1 = Q(y) = 0$. By replacing $x$ by $y$, we may assume that $\det x_1 = Q(x) = 0$.

Suppose $x_1 = 0$. Then if $\det x_2 \neq 0$, $x \in Y_{2k}^{\mathrm{ss}}$. Since $x \neq 0$, $x_2 \neq 0$. So if $\det x_2 = 0$, the rank of $x_2$ is one. This implies that there exists $g_1 \in G_{1k}$ such that $g_1 x_2 {}^t g_1^{\sigma}$ is of the form $\begin{pmatrix} 0 & 0 \\ 0 & c \end{pmatrix}$ where $c \neq 0$. Then $x \in S_{3k}$.

Suppose $x_1 \neq 0$. Then the rank of $x_1$ is one. By changing $x$ if necessary, we may assume that $x_{10} = x_{11} = 0$, $x_{12} \neq 0$. Since $Q(x) = 0$, $x_{20} = 0$. If $x_{21} = 0$, by changing $x$ by an element of $G_{2k}$ if necessary, we may assume $x_{12} = 0$. Then $x \in Z_{3k}^{\mathrm{ss}}$. If $x_{21} \neq 0$, $x \in Y_{1k}^{\mathrm{ss}}$. This proves $V_k \setminus \{0\} = V_k^{\mathrm{ss}} \cup S_{1k} \cup S_{2k} \cup S_{3k}$.

If $x \in S_{2k}$ or $S_{3k}$, the dimension of the subspace of $W$ spanned by $x_1, x_2$ is one. If $x \in S_{1k}$, the dimension of the subspace of $W$ spanned by $x_1, x_2$ is two. So $S_{1k} \cap S_{2k} = S_{1k} \cap S_{3k} = \emptyset$. If $x \in S_{2k}$, $\det(v_1 x_1 + v_2 x_2) \neq 0$, but if $x \in S_{3k}$, $\det(v_1 x_1 + v_2 x_2) = 0$. So $S_{2k} \cap S_{3k} = \emptyset$. This proves (1).

It is easy to see that $P_{ik}Y_{ik}^{\mathrm{ss}} = Y_{ik}$ for all $i$. So in order to prove (2), it is enough to show that if $x \in Y_{ik}^{\mathrm{ss}}, g = (g_1, g_2) \in G_k$ and $gx \in Y_{ik}^{\mathrm{ss}}$, then $g \in P_{ik}$.

Consider the case $i = 2$. By the Bruhat decomposition of $G_2$, $g_2 \in B_{2k}$ or $B_{2k}\nu B_{2k}$. Suppose $g_2 \in B_{2k}\nu B_{2k}$. If $b_2 \in B_{2k}$, $b_2 x \in Y_{2k}^{\mathrm{ss}}$. So if $y = (y_1, y_2) = gx$, $y_{1j} \neq 0$ for a certain $j$, which is not possible. So $g \in P_{1k}$.

If $i = 1, 3$, $g \in B_k$, $B_k(\nu, 1)B_k$, $B_k(1, \nu)B_k$, or $B_k(\nu, \nu)B_k$. Let $y = (y_1, y_2) = gx$. Consider the case $i = 1$. If $g \in B_k(\nu, 1)B_k$, $y_{10} \neq 0$. If $g \in B_k(1, \nu)B_k$, $y_{11} \neq 0$. Suppose $g = b(\nu, \nu)b'$ where $b, b' \in B_k$. Let $z = (z_1, z_2) = (\nu, \nu)b'x$. Then $z_{11} \neq 0$. If $z_{10} \neq 0$, $y_{10} \neq 0$ also. If $z_{10} = 0$, $y_{11} \neq 0$. So $y \notin Y_{2k}^{\mathrm{ss}}$, which is a contradiction. This implies $g \in B_k$.

Consider the case $i = 3$. If $g \in B_k(\nu, 1)B_k$, $y_{20} \neq 0$. If $g \in B_k(1, \nu)B_k$, $y_{12} \neq 0$. If $g \in B_k(\nu, \nu)B_k$, $y_{10} \neq 0$. So $g \in B_k$. This proves (2).

Q.E.D.

Let $\widetilde{T} = \mathrm{Ker}(G \to \mathrm{GL}(V))$. For $x \in V^{\mathrm{ss}}$, let $G_x$ be the stabilizer of $x$, and $G_x^0$ its connected component of 1. We define

$$L_0 = \{x \in V_k^{\mathrm{ss}} \mid G_x^0/\widetilde{T} \text{ does not have a non-trivial rational character }\}.$$

Note that we proved in [1] that the set $G_k \setminus V_k^{\mathrm{ss}}$ corresponds bijectively with fields $k'/k$ such that $[k' : k] \leq 2$. Moreover, if $x \in V_k^{\mathrm{ss}}$, the corresponding field $k'$ is generated over $k$ by the roots of $F_x(v)$. We proved in §1 of [1] that $L_0$ is the set of points such that $F_x(v)$ does not factor over $k_1$. So $V_k^{\mathrm{ss}} \setminus L_0$ consists of precisely two $G_k$–orbits. They are represented by the following two elements

$$w_1 = \left( \begin{pmatrix} 0 & 0 \\ 0 & 1 \end{pmatrix}, \begin{pmatrix} 1 & 0 \\ 0 & 0 \end{pmatrix} \right),$$
$$w_2 = \left( \begin{pmatrix} 0 & 1 \\ 1 & 0 \end{pmatrix}, \begin{pmatrix} 0 & \alpha_0 \\ -\alpha_0 & 0 \end{pmatrix} \right).$$

Note that

$$F_{w_1}(v) = v_1 v_2,\ F_{w_2}(v) = v_1^2 - \alpha_0^2 v_2^2,$$

which implies that $w_1$ corresponds to the field $k$ and $w_2$ corresponds to the field $k_1$. We define $V_{\mathrm{st},ik}^{\mathrm{ss}} = G_k w_i$ for $i = 1, 2$. Then

$$V_k^{\mathrm{ss}} = L_0 \coprod V_{\mathrm{st},1k}^{\mathrm{ss}} \coprod V_{\mathrm{st},2k}^{\mathrm{ss}}.$$

We determine the structure of $V_{\mathrm{st},i}^{\mathrm{ss}}$ for $i = 1, 2$. Let $H_1 \subset G_k$ be the subgroup generated by $T_k$ and $(\nu, \nu)$, and $H_2 \subset G_k$ the subgroup generated by $T_{1k} \times G_{2k}$ and $(\nu, 1)$. Let

$$\begin{aligned}
Y_{0,1} &= \{x \in V \mid x_{10} = x_{11} = 0\}, \\
Z_{0,1} &= \{x \in V \mid x_{10} = x_{11} = x_{22} = 0\}, \\
Z'_{0,1} &= \{x \in V \mid x_{10} = x_{11} = x_{21} = x_{22} = 0\}, \\
{Z'}_{0,1}^{\mathrm{ss}} &= \{x \in Z'_{0,1} \mid x_{12}, x_{20} \neq 0\}, \\
Z_{0,2} &= \{x \in V \mid x_{10} = x_{20} = 0\}, \\
Z'_{0,2} &= \{x \in V \mid x_{10} = x_{20} = x_{12} = x_{22} = 0\}, \\
{Z'}_{0,2}^{\mathrm{ss}} &= \{x \in Z'_{0,2} \mid \{x_{11}, x_{21}\} \text{ is linearly independent over } k\}.
\end{aligned}$$

**Proposition (1.2)** $V_{\mathrm{st},ik}^{\mathrm{ss}} \cong G_k \times_{H_{ik}} {Z'}_{0,ik}^{\mathrm{ss}}$.

*Proof.* The proof of the case $i=1$ is easy and is left to the reader. Consider the case $i=2$. Clearly, $w_2 \in {Z'}_{0,ik}^{\mathrm{ss}}$ and $H_{2k}$ stabilizes the set ${Z'}_{0,2k}^{\mathrm{ss}}$. So it is enough to show that if $x, y \in {Z'}_{0,ik}^{\mathrm{ss}}$, $g \in G_k$ and $gx = y$, then $g \in H_{2k}$.

Note that $G_{k_1} \cong \mathrm{GL}(2)_{k_1} \times \mathrm{GL}(2)_{k_1} \times \mathrm{GL}(2)_{k_1}$ and the action of the non-trivial element $\sigma \in \mathrm{Gal}(k_1/k)$ is given by $(g_1, g_2, g_3)^\sigma = (g_2^\sigma, g_1^\sigma, g_3^\sigma)$. By this identification, the inclusion $G_k \to G_{k_1}$ can be identified with the map $(g_1, g_2) \to (g_1, g_1^\sigma, g_2)$. Also $W \otimes k_1 \cong \mathrm{M}(2,2)_{k_1}$ and the action of $G_{k_1}$ on $V \otimes k_1$ can be considered as

$$(g_1, g_2, g_3) M_x(v) = g_1 M_x(v g_3) {}^t g_2$$

for $M_x(v) = v_1 x_1 + v_2 x_2 \in \mathrm{M}(2,2)_{k_1} \otimes k^2$.

Let

$$g_{w_2} = \left( \nu, I_2, \begin{pmatrix} 1 & 1 \\ \alpha_0 & -\alpha_0 \end{pmatrix} \right).$$

Then $g_{w_2} w_1 = w_2$. This implies $G^0_{w_2 k_1} = g_{w_2} G^0_{w_1 k_1} g_{w_2}^{-1}$. It is easy to see that

$$G^0_{w_1 k_1} = \left\{ (a(t_{11}, t_{12}), a(t_{21}, t_{22}), a(t_{31}, t_{32})) \;\middle|\; \begin{array}{c} t_{ij} \in k_1^\times \text{ for all } i, j, \\ t_{11} t_{21} t_{31} = t_{12} t_{22} t_{32} = 1 \end{array} \right\}.$$

Let $t = (a(t_{11}, t_{12}), a(t_{21}, t_{22}), a(t_{31}, t_{32})) \in G^0_{w_1 k_1}$. Then the element $g_{w_2} t g_{w_2}^{-1}$ belongs to $G^0_{wk}$ if and only if $(g_{w_2} t g_{w_2}^{-1})^\sigma = g_{w_2} t g_{w_2}^{-1}$. This condition is equivalent to

$$((g_{w_2}^\sigma)^{-1} g_{w_2})^{-1} t^\sigma (g_{w_2}^\sigma)^{-1} g_{w_2} = t. \tag{1.3}$$

It is easy to see that

$$(g_{w_2}^\sigma)^{-1} = \left( I_2, \nu, \frac{1}{2\alpha_0} \begin{pmatrix} \alpha_0 & -1 \\ \alpha_0 & 1 \end{pmatrix} \right),$$
$$(g_{w_2}^\sigma)^{-1} g_{w_2} = (\nu, \nu, \nu).$$

Since

$$(\nu, \nu, \nu) t^\sigma (\nu, \nu, \nu) = (a(t_{22}^\sigma, t_{21}^\sigma), a(t_{12}^\sigma, t_{11}^\sigma), a(t_{32}^\sigma, t_{31}^\sigma)),$$

(1.3) is satisfied if and only if $t_{21} = t_{12}^\sigma$, $t_{22} = t_{11}^\sigma$, $t_{32} = t_{31}^\sigma$. Therefore,

$$G^0_{w_2 k} = g_{w_2} \{ (a(t_{11}, t_{12}), a(t_{12}^\sigma, t_{11}^\sigma), a(t_{11}^{-1} (t_{12}^\sigma)^{-1}, (t_{11}^\sigma)^{-1} t_{12}^{-1})) \mid t_{11}, t_{12} \in k_1^\times \} g_{w_2}^{-1}.$$

Let

$$A(t_{11}, t_{12}) = \begin{pmatrix} 1 & 1 \\ \alpha_0 & -\alpha_0 \end{pmatrix} a(t_{11}^{-1} (t_{12}^\sigma)^{-1}, (t_{11}^\sigma)^{-1} t_{12}^{-1}) \begin{pmatrix} 1 & 1 \\ \alpha_0 & -\alpha_0 \end{pmatrix}^{-1}.$$

Then $A(t_{11}, t_{12}) \in \mathrm{GL}(2)_k$. Since

$$g_{w_2} (a(t_{11}, t_{12}), a(t_{12}^\sigma, t_{11}^\sigma), 1) g_{w_2}^{-1} = (a_2(t_{12}, t_{11}), a_2(t_{12}, t_{11})^\sigma, 1),$$

we get

$$G^0_{w_2 k} = \{(a_2(t_{12}, t_{11}), A(t_{11}, t_{12})) \mid t_{11}, t_{12} \in k_1^\times\}.$$

Note we are identifying $(a_2(t_{12}, t_{11}), a_2(t_{12}, t_{11})^\sigma)$ with $a_2(t_{12}, t_{11})$.

Since

$$g_{w_2}(\nu, \nu, \nu) g_{w_2}^{-1} = (\nu, \nu, a(1, -1)),$$

(which corresponds to $(\nu, a(1, -1))$), $G_{w_2 k}$ is generated by $G^0_{w_2 k}$ and $(\nu, a(1, -1))$.

Suppose $gx = y$ where $x, y \in Z'^{\rm ss}_{0,2k}$, $g \in G_k$. It is easy to see that there exist $h_x, h_y \in G_{2k}$ such that $x = h_x w_2$, $y = h_y w_2$. So $g h_x w_2 = h_y w_2$. This implies $h_y^{-1} g h_x \in G_{w_2 k} \subset H_{2k}$. Since $h_x, h_y \in H_{2k}$, this completes the proof of Proposition (1.2).

Q.E.D.

## §2 The zeta function and the smoothed Eisenstein series

In this section we define the zeta function and discuss basic properties of the smoothed Eisenstein series.

Let

$$(2.1) \qquad a(t_{21}, t_{22}) = \begin{pmatrix} t_{21} & 0 \\ 0 & t_{22} \end{pmatrix}, \; n(u) = \begin{pmatrix} 1 & 0 \\ u & 1 \end{pmatrix}.$$

where $t_{21}, t_{22} \in \mathbb{A}_{k_1}^\times$, $u \in \mathbb{A}_{k_1}$. Throughout this paper, we use the following notation.

$$(2.2) \qquad \begin{aligned} \widehat{t}^0 &= (a_2(t_{11}, t_{12}), a_2(t_{21}, t_{22})), \\ d(\lambda_1, \lambda_2) &= (a_2(\underline{\lambda}_{1k_1}^{-1}, \underline{\lambda}_{1k_1}), a_2(\underline{\lambda}_2^{-1}, \underline{\lambda}_2)), \; t^0 = d(\lambda_1, \lambda_2)\widehat{t}^0, \\ n(u) &= (n(u_1), n(u_2)) \end{aligned}$$

for

$$t_{11}, t_{12} \in \mathbb{A}^1_{k_1}, \; t_{21}, t_{22} \in \mathbb{A}^1, \; \lambda_1, \lambda_2 \in \mathbb{R}_+, \; u = (u_1, u_2), \; u_1 \in \mathbb{A}_{k_1}, u_2 \in \mathbb{A}.$$

Let

$$(2.3) \qquad \begin{aligned} G^0_{1\mathbb{A}} &= \{g_1 \in G_{1\mathbb{A}} \mid |\det g_1|_{k_1} = 1\}, \\ G^0_{2\mathbb{A}} &= \{g_2 \in G_{2\mathbb{A}} \mid |\det g_2| = 1\}, \\ G^0_{\mathbb{A}} &= G^0_{1\mathbb{A}} \times G^0_{2\mathbb{A}}, \; \widetilde{G}_{\mathbb{A}} = \mathbb{R}_+ \times G^0_{\mathbb{A}}, \\ \widehat{T}^0_{1\mathbb{A}} &= \{a(t_{11}, t_{12}) \mid t_{11}, t_{12} \in \mathbb{A}^1_{k_1}\}, \\ \widehat{T}^0_{2\mathbb{A}} &= \{a(t_{21}, t_{22}) \mid t_{21}, t_{22} \in \mathbb{A}^1\}, \\ T^0_{1\mathbb{A}} &= \{a(\underline{\lambda}_{1k_1}^{-1}, \underline{\lambda}_{1k_1}) t_1 \mid \lambda_1 \in \mathbb{R}_+, \; t_1 \in \widehat{T}^0_{1\mathbb{A}}\}, \\ T^0_{2\mathbb{A}} &= \{a(\underline{\lambda}_2^{-1}, \underline{\lambda}_2) t_2 \mid \lambda_2 \in \mathbb{R}_+, \; t_2 \in \widehat{T}^0_{2\mathbb{A}}\}, \\ \widehat{T}^0_{\mathbb{A}} &= \{\widehat{t}^0 \mid t_{11}, t_{12} \in \mathbb{A}^1_{k_1}, \; t_{21}, t_{22} \in \mathbb{A}^1\} = \widehat{T}^0_{1\mathbb{A}} \times \widehat{T}^0_{2\mathbb{A}}, \\ T^0_{\mathbb{A}} &= \{d(\lambda_1, \lambda_2)\widehat{t}^0 \mid \lambda_1, \lambda_2 \in \mathbb{R}_+, \; \widehat{t}^0 \in \widehat{T}^0_{\mathbb{A}}\} = T^0_{1\mathbb{A}} \times T^0_{2\mathbb{A}}. \end{aligned}$$

The group $\widetilde{G}_{\mathbb{A}}$ acts on $V_{\mathbb{A}}$ by assuming that $\lambda \in \mathbb{R}_+$ acts by multiplication by $\underline{\lambda}$.

If $t^0 \in T^0_{\mathbb{A}}$ and $x = (x_{ij}) \in V_k$, there are elements $\gamma_{ij}(t^0) \in \mathbb{A}^\times$ for $i = 1, 2$, $j = 0, 2$ and $\gamma_{i1}(t^0) \in \mathbb{A}^\times_{k_1}$ for $i = 1, 2$ such that $t^0 x = (\gamma_{ij}(t^0) x_{ij})$. It is easy to see that

$$\begin{aligned}
&\gamma_{10}(d(\lambda_1,\lambda_2)) = \underline{\lambda}_1^{-1}\underline{\lambda}_2^{-1},\ \gamma_{11}(d(\lambda_1,\lambda_2)) = \underline{\lambda}_2^{-1},\ \gamma_{12}(d(\lambda_1,\lambda_2)) = \underline{\lambda}_1\underline{\lambda}_2^{-1},\\
&\gamma_{20}(d(\lambda_1,\lambda_2)) = \underline{\lambda}_1^{-1}\underline{\lambda}_2,\ \gamma_{21}(d(\lambda_1,\lambda_2)) = \underline{\lambda}_2,\ \gamma_{22}(d(\lambda_1,\lambda_2)) = \underline{\lambda}_1\underline{\lambda}_2.
\end{aligned}$$

Throughout this paper, we express elements $\widetilde{g} \in \widetilde{G}_{\mathbb{A}}$, $g^0 \in G^0_{\mathbb{A}}$ as

$$\widetilde{g} = (\lambda, g_1, g_2),\ g^0 = (g_1, g_2), \tag{2.4}$$

where $\lambda \in \mathbb{R}_+$, $g_1 \in G^0_{1\mathbb{A}}$, and $g_2 \in G^0_{2\mathbb{A}}$. We identify the element $g^0 \in G^0_{\mathbb{A}}$ with $(1, g^0)$ and $g_1 \in G_{1\mathbb{A}}, g_2 \in G^0_{2\mathbb{A}}$ with $(1, g_1, 1), (1, 1, g_2)$. Let $K_i \subset G_{i\mathbb{A}}$ be the standard maximal compact subgroup for $i = 1, 2$, and $K = K_1 \times K_2$. Then $K$ is a maximal compact subgroup of $G_{\mathbb{A}}$. Let $dk$ be the Haar measure on $K$ such that $\int_K dk = 1$. Let

$$b^1 = d(\lambda_1, \lambda_2)\widehat{t}^0 n(u), \tag{2.5}$$

where $d(\lambda_1, \lambda_2)$, etc. are as in (2.2). Let $B^1_{\mathbb{A}}$ be the subgroup of $G_{\mathbb{A}}$ consisting of such $b^1$'s.

We use the Haar measure $d^\times t'$ on $t' \in \mathbb{A}^1$ or $\mathbb{A}^1_{k_1}$ such that the volume of $\mathbb{A}^1/k^\times$ or $\mathbb{A}^1_{k_1}/k_1^\times$ is 1. If $t = \underline{\lambda} t'$ or $\underline{\lambda}_{k_1} t'$ where $t' \in \mathbb{A}^1$ or $\mathbb{A}^1_{k_1}$, we use $d^\times t = d^\times \lambda d^\times t'$ as the Haar measure on $\mathbb{A}^\times$ or $\mathbb{A}^\times_{k_1}$. We do not use the subscript $k_1$ for the measure on $\mathbb{A}^1_{k_1}$ because the situation will be clear from the context. Let $du_1$, $du_2$ be Haar measures on $\mathbb{A}_{k_1}$, $\mathbb{A}$ respectively such that the volumes of $\mathbb{A}_{k_1}/k_1$, $\mathbb{A}/k$ are 1. We define

$$\begin{aligned}
&d^\times \widehat{t}^0 = d^\times t_{11} d^\times t_{12} d^\times t_{21} d^\times t_{22},\ d^\times t^0 = d^\times \lambda_1 d^\times \lambda_2 d^\times \widehat{t}^0,\\
&du = du_1 du_2,\ db^1 = \lambda_1^2 \lambda_2^2 d^\times t^0 du.
\end{aligned} \tag{2.6}$$

We use $dg^0 = dk db^1$ as the Haar measure on $G^0_{\mathbb{A}}$. Let $dg_1$, $dg_2$ be the Haar measures on $G^0_{1\mathbb{A}}$, $G^0_{2\mathbb{A}}$ which we defined in §1.1 [5] (replacing $k$ by $k_1$ for $dg_1$). Then $dg^0 = dg_1 dg_2$. We define $d\widetilde{g} = d^\times \lambda dg^0$.

For $\eta > 0$, we define

$$\begin{aligned}
&T^0_{1\eta+} = \{d(\lambda_1, 1) \mid \lambda_1 \in \mathbb{R}_+,\ \lambda_1 \le \eta\},\\
&T^0_{2\eta+} = \{d(1, \lambda_2) \mid \lambda_2 \in \mathbb{R}_+,\ \lambda_2 \le \eta\},\\
&T^0_{\eta+} = \{d(\lambda_1, \lambda_2) \mid \lambda_1, \lambda_2 \in \mathbb{R}_+,\ \lambda_1, \lambda_2 \le \eta\} = T^0_{1\eta+} \times T^0_{2\eta+}.
\end{aligned} \tag{2.7}$$

Let $\Omega_i \subset \widehat{T}^0_{i\mathbb{A}} N_{i\mathbb{A}}$ be a compact subset for $i = 1, 2$. We define $\mathfrak{S}^0_i = K_i T^0_{i\eta+} \Omega_i$ for $i = 1, 2$ and $\Omega = \Omega_1 \times \Omega_2$. The sets $\mathfrak{S}^0_1, \mathfrak{S}^0_2, \mathfrak{S}^0$ are called Siegel sets. It is known that for a suitable choice of $\eta$ and $\Omega_i$, $\mathfrak{S}^0_i$ surjects to $G^0_{i\mathbb{A}}/G_{ik}$ for $i = 1, 2$. Also there exists another compact set $\widehat{\Omega}_i \subset G^0_{i\mathbb{A}}$ such that $\mathfrak{S}^0_i \subset \widehat{\Omega}_i T^0_{i\eta+}$ for $i = 1, 2$. Let $\widehat{\Omega} = \widehat{\Omega}_1 \times \widehat{\Omega}_2$.

For $x=(x_1,x_2)=(x_{ij})$ and $y=(y_1,y_2)=(y_{ij})$, we define

$$\begin{aligned}[x,y]' &= \mathrm{Tr}_{k_1/k}(\mathrm{tr}(x_1y_1))+\mathrm{Tr}_{k_1/k}(\mathrm{tr}(x_2y_2))\\ &= x_{10}y_{10}+\mathrm{Tr}_{k_1/k}(x_{11}y_{11}^{\sigma})+x_{12}y_{12}+x_{20}y_{20}+\mathrm{Tr}_{k_1/k}(x_{21}y_{21}^{\sigma})+x_{22}y_{22}.\end{aligned}$$

This is a non-degenerate bilinear form on $V$ defined over $k$. By an easy consideration, $[gx,y]'=[x,{}^tgy]'$ for all $g,x,y$. We define $[x,y]=[x,(\nu,\nu)y]'$. For $\widetilde{g}=(\lambda,g_1,g_2)\in\widetilde{G}_{\mathbb{A}}$, we define

$$\widetilde{g}^{\iota}=(\lambda^{-1},\nu{}^tg_1^{-1}\nu,\nu{}^tg_2^{-1}\nu).$$

This is an involution and the above bilinear form satisfies $[\widetilde{g}x,y]=[x,(\widetilde{g}^{\iota})^{-1}y]$.

For $\Phi\in\mathscr{S}(V_{\mathbb{A}})$, we define its Fourier transform by

$$\widehat{\Phi}(x)=\int_{V_{\mathbb{A}}}\Phi(y)<[x,y]>dy.$$

It is easy to see that the Fourier transform of $\Phi(\widetilde{g}\cdot\quad)$ is $\lambda^{-8}\widehat{\Phi}(g^{\iota}\cdot\quad)$.

**Definition (2.8)** *For any $G_k$–invariant subset $\subset V_k$ and $\Phi\in\mathscr{S}(V_{\mathbb{A}})$, we define*

$$\Theta_S(\Phi,\widetilde{g})=\sum_{x\in S}\Phi(\widetilde{g}x).$$

We are going to define the zeta function for the prehomogeneous vector space in this paper now. The zeta function for our case is not defined for $V_k^{\mathrm{ss}}$ and we have to consider $L_0$ instead of $V_k^{\mathrm{ss}}$.

**Lemma (2.9)** (1) *There exists a slowly increasing function $h(\lambda,\lambda_1,\lambda_2)$ such that for any $N_1,N_2\geq 1$,*

$$\Theta_{L_0}(\Phi,\widetilde{g})\ll h(\lambda,\lambda_1,\lambda_2)\sup((\lambda\lambda_1^{-1}\lambda_2^{-1})^{-N_1},(\lambda\lambda_2^{-1})^{-2N_1}(\lambda\lambda_1^{-1}\lambda_2)^{-N_2})$$

*for $\widetilde{g}\in\mathbb{R}_+\times\mathfrak{G}^0$.*

*Proof.* We first prove that if $x\in L_0$, either $x_{10}\neq 0$, or $x_{10}=0,x_{11},x_{20}\neq 0$.

Suppose $x_{10}=0$. If $x_{11}=0$, $\det x_1=0$, which means $F_x(v)$ factors over $k\subset k_1$. So we may assume that $x_{11}\neq 0$. Suppose $x_{20}=0$. By applying an element of $G_{k_1}$ (not $G_k$), we can make $x_{21}=0$. Then $\det x_2=0$ and $F_x(v)$ factors. Whether or not $F_x(v)$ factors over $k_1$ does not change by applying an element of $G_{k_1}$ to $x$. So this means if $x_{10}=x_{20}=0$, $F_x(v)$ factors over $k_1$. This proves that if $x\in L_0$, either $x_{10}\neq 0$, or $x_{10}=0,x_{11},x_{20}\neq 0$.

Now the statement of the lemma is a direct consequence of Lemmas (1.2.3), (1.2.8) [5]. Note that the coordinate $x_{21}$ belongs to $k_1$ which is a two dimensional vector space over $k$ and this is why we get $(\lambda\lambda_2^{-1})^{-2N_1}$ instead of $(\lambda\lambda_2^{-1})^{-N_1}$.

Q.E.D.

Let $\omega_1,\omega_2$ be characters of $\mathbb{A}_{k_1}^{\times}/k_k^{\times}$, $\mathbb{A}^{\times}/k^{\times}$ respectively. We put $\omega=(\omega_1,\omega_2)$. For this $\omega$, we define a character of $\widetilde{G}_{\mathbb{A}_k}/G_k$ by $\omega(\widetilde{g})=\omega_1(\det g_1)\omega_2(\det g_2)$. We define $\delta(\omega_1) = 1$ if $\omega_1$ is a trivial character and $\delta(\omega_1) = 0$ otherwise. We define $\delta(\omega_2)$ similarly. Also we define $\delta_{\#}(\omega) = \delta(\omega_1)\delta(\omega_2)$.

**Definition (2.10)** *For $\Phi \in \mathscr{S}(V_{\mathbb{A}})$, $\omega$ as above and a complex variable $s$, we define*

$$(1) \qquad Z(\Phi, \omega, s) = \int_{\widetilde{G}_{\mathbb{A}_k}/G_k} \lambda^s \omega(\widetilde{g}) \Theta_{L_0}(\Phi, \widetilde{g}) d\widetilde{g},$$

$$(2) \qquad Z_+(\Phi, \omega, s) = \int_{\substack{\widetilde{G}_{\mathbb{A}_k}/G_k \\ \lambda \geq 1}} \lambda^s \omega(\widetilde{g}) \Theta_{L_0}(\Phi, \widetilde{g}) d\widetilde{g}.$$

By Lemma (2.9), the integral (1) converges absolutely and locally uniformly on a certain right half plane and the integral (2) is an entire function.

For $\lambda \in \mathbb{R}_+$ and $\Phi \in \mathscr{S}(V_{\mathbb{A}})$, we define $\Phi_\lambda(x) = \Phi(\underline{\lambda} x)$. Let

$$(2.11) \qquad \begin{aligned} J(\Phi, g^0) &= \sum_{x \in V_k \setminus L_0} \widehat{\Phi}((g^0)^\iota x) - \sum_{x \in V_k \setminus L_0} \Phi(g^0 x), \\ I^0(\Phi, \omega) &= \int_{G^0_{\mathbb{A}}/G_k} \omega(g^0) J(\Phi, g^0) dg^0, \\ I(\Phi, \omega, s) &= \int_0^1 \lambda^s I^0(\Phi_\lambda, \omega) d^\times \lambda. \end{aligned}$$

Then by the Poisson summation formula,

$$Z(\Phi, \omega, s) = Z_+(\Phi, \omega, s) + Z(\widehat{\Phi}, \omega^{-1}, 8 - s) + I(\Phi, \omega, s).$$

We define $M_\omega \Phi$ as in (3.1.11) [5]. Since $Z(\Phi, \omega, s) = Z(M_\omega \Phi, \omega, s)$, we assume that $\Phi = M_\omega \Phi$ for the rest of this paper.

For the rest of this section, we discuss the smoothed Eisenstein series. Let $z_1 = (z_{11}, z_{12}), z_2 = (z_{21}, z_{22}) \in \mathbb{C}^2$ where $z_{11} + z_{12} = z_{21} + z_{22} = 0$, and $z = (z_1, z_2)$. For $t_{11}, t_{12} \in \mathbb{A}^\times_{k_1}$, $t_{21}, t_{22} \in \mathbb{A}^\times$, we define

$$(a(t_{11}, t_{12}), a(t_{21}, t_{22}))^z = |t_{11}|_{k_1}^{z_{11}} |t_{12}|_{k_1}^{z_{12}} |t_{21}|^{z_{21}} |t_{22}|^{z_{22}}.$$

Let $E_{B_1}(g_1, z_1), E_{B_2}(g_2, z_2)$ be the Eisenstein series of $G_{1\mathbb{A}}, G_{2\mathbb{A}}$ with respect to $B_1, B_2$ respectively. Then $E_B(g^0, z) = E_{B_1}(g_1, z_1) E_{B_2}(g_2, z_2)$ is the Eisenstein series of $G$ with respect to $B$. Let $\rho_1 = \rho_2 = (\frac{1}{2}, -\frac{1}{2})$ and $\rho = (\rho_1, \rho_2)$. We can consider $\rho$ as half the sum of positive weights. Let $\psi(z)$ be an entire function which is rapidly decreasing with respect to $\mathrm{Im}(z)$ on any vertical strip. Moreover, we assume that

$$(2.12) \qquad \psi(z_{12}, z_{11}, z_{22}, z_{21}) = \psi(z_{11}, z_{12}, z_{21}, z_{22}), \ \psi(\rho) \neq 0.$$

We choose a constant $C > 4$. Let

$$L(z) = (z_{11} - z_{12}) + C(z_{21} - z_{22}), \ \Lambda(w; z) = \frac{\psi(z)}{w - L(z)},$$

$$dz_1 = d(z_{11} - z_{12}), \ dz_2 = d(z_{21} - z_{22}), \ dz = dz_1 dz_2.$$

We define the smoothed Eisenstein series as follows.

$$\mathscr{E}(g^0, w) = \left(\frac{1}{2\pi\sqrt{-1}}\right)^2 \int_{\mathrm{Re}(z)=q} E_B(g^0, z)\Lambda(w; z)dz, \tag{2.13}$$

where $q = (q_{11}, q_{12}, q_{21}, q_{22}) \in \mathbb{R}^4$ satisfies $q_{11} - q_{12} > 1$, $q_{21} - q_{22} > 1$.

In [5], we developed a general procedure to compute the principal part of the zeta function using $\mathscr{E}(g^0, w)$ assuming that the group is a product of $\mathrm{GL}(n)$'s over the ground field $k$, Even though our group is not of this form, the Eisenstein series (without smoothing) we are using are those of $\mathrm{GL}(2)_{k_1}$, $\mathrm{GL}(2)_k$, and the estimates of these Eisenstein series needed to prove statements in Chapter 3 [5] are covered in Chapter 2 [5] (by replacing $k$ by $k_1$). The only issue we have to worry about is that the convex hull construction of the stratification works in every step. Since this is true for our case, we will use results in Chapter 3 [5]. Here we summarize for which situation results in Chapter 3 [5] hold.

**Assumptions to apply results of Chapter 3 [5]**

Suppose $G$ is a product of $\mathrm{GL}(n)$'s over finite extensions over $k$. We define $T, N, B \subset G$ similarly as in Chapters 2,3 [5] ($T$ is the maximal torus, $B$ is the Borel subgroup, and $N$ is its unipotent radical). Let $\overline{T} \subset T$ be the maximal split torus. For example, if $G = \mathrm{GL}(n)_{k_1}$ where $k_1/k$ is a finite extension, we can choose $\overline{T}$ to be the set of diagonal matrices with entries in $k^\times$. Let $(G, V, \chi)$ be a prehomogeneous vector space. Let $T' = \overline{T} \cap \mathrm{Ker}(\chi)$, $\mathfrak{t} = X_*(T') \otimes \mathbb{R}$, and $\mathfrak{t}^* = X^*(T') \otimes \mathbb{R}$, where $X_*(T'), X^*(T')$ are the groups of rational one parameter subgroups, rational characters respectively. We consider an inner product on $\mathfrak{t}^*$ which is invariant under the action of the Weyl group. By this metric, we can identify $\mathfrak{t}$ with $\mathfrak{t}^*$ and we can consider any $\beta \in \mathfrak{B}$ as an element of $\mathfrak{t}$ also. Since $\overline{T}$ is split and number fields are characteristic zero fields, we can still diagonalize the action of $\overline{T}$. So we choose a coordinate system $x = (x_i)$ such that $t \cdot x = (\gamma_i(t)x_i)$ for $t \in \overline{T}$ where $\gamma_i$ is a rational character of $\overline{T}$ for all $i$. By restricting to $T'$, we can consider $\gamma_i \in \mathfrak{t}^*$. We define the set $\mathfrak{B}$ of minimal combination of weights using Definition (3.2.3) [5]. For each $\beta \in \mathfrak{B}$, we can construct $M_\beta, P_\beta, S_\beta, Y_\beta^{\mathrm{ss}}, Z_\beta^{\mathrm{ss}}$ in exactly the same manner as in §3.2 [5]. Note that the notion of semi-stable points does not change by field extensions. Also note that $M_\beta$ is again a product of $\mathrm{GL}(n)$'s over finite extensions of $k$. So we can consider the notion of $\beta$-sequences using Definition (3.2.6) [5] (note that there is a typo in the statement of (3.2.6)(2) [5] and $\cap_{1\le i\le j} M_{\beta_i}$ must be replaced by $\cap_{1\le i<j} M_{\beta_i}$). We can also construct $M_{\mathfrak{d}}, P_{\mathfrak{d}}, S_{\mathfrak{d}}, Y_{\mathfrak{d}}^{\mathrm{ss}}, Z_{\mathfrak{d}}^{\mathrm{ss}}$ as in §3.2 [5].

Suppose $\mathfrak{d} = (\beta_1, \cdots, \beta_a)$ is a $\beta$-sequence. Consider $\beta$-sequences of the form $\mathfrak{d}' = (\beta_1, \cdots, \beta_a, \beta_{a+1})$.

**Assumption (2.14)** (1) $Z_{\mathfrak{d}k} \setminus \{0\} = Z_{\mathfrak{d}k}^{\mathrm{ss}} \coprod_{\beta_{a+1}} S_{\mathfrak{d}'k}$.
(2) $S_{\mathfrak{d}'k} \cong M_{\mathfrak{d}k} \times_{P_{\mathfrak{d}'k}} Y_{\mathfrak{d}'k}^{\mathrm{ss}}$ for all $\beta_{a+1}$.

We assumed in [5] that we only get prehomogeneous vector spaces in induction steps. So we make this assumption also.

**Assumption (2.15)** For all the $\beta$-sequences $\mathfrak{d}$, $(M_{\mathfrak{d}}, Z_{\mathfrak{d}})$ is a prehomogeneous vector space.

Assumption (2.14) is probably always satisfied, but we do not consider this issue in this paper.

Suppose Assumptions (2.14), (2.15) are satisfied. We define $\overline{T}_+$ in the usual manner. If $k_1/k$ is a field extension of degree $n$ and $\lambda \in \mathbb{R}_+$, we define $\underline{\lambda}_{k_1} \in \mathbb{A}_{k_1}^\times$, $\underline{\lambda} \in \mathbb{A}^\times$ in the same manner as before. Then $\underline{\lambda} = \underline{\lambda}_{k_1}^n$. We use $\overline{T}_+$ in place of $T_+$ in [5]. If $\mathfrak{d} = (\beta_1, \cdots, \beta_a)$ is a $\beta$-sequence, we can define $\nu_{\mathfrak{d}1}, \cdots, \nu_{\mathfrak{d}a}$, $M_{\mathfrak{d}\mathbb{A}}^0$, etc. in the same way as in [5, p. 72]. Also the sets $A_{\mathfrak{p}0}$, etc. can be defined without any modification. Therefore, all the results in Chapters 2, 3 [5] hold.

We go back to our case. We proved a stratification of $V_k$ in Proposition (1.1). It is easy to see that this is precisely the one given by the convex hull construction in §3.2 [5] using the convex hulls in [5, p. 153] (with the horizontal size changed). Let $\beta_1, \beta_2, \beta_3 \in \mathfrak{t}^*$ be the corresponding points. In the stratification of $Z_{\beta_2 k}$, we get a stratum which corresponds to $\beta_3$. So Assumptions (2.14), (2.15) are satisfied for our case. Let $\mathfrak{d}_1 = (\beta_1), \mathfrak{d}_2 = (\beta_2), \mathfrak{d}_3 = (\beta_3), \mathfrak{d}_4 = (\beta_1, \beta_3)$. These are all the $\beta$-sequences we have to consider.

We recall properties of $\mathscr{E}(g^0, w)$. Let $s_1 = z_{11} - z_{12}$, $s_2 = z_{21} - z_{22}$. We define

$$
\begin{aligned}
(2.16)\qquad & \phi_{k_1}(z_1) = \phi_{k_1}(s_1) = Z_{k_1}(s_1) Z_{k_1}(s_1+1)^{-1}, \\
& \phi(z_2) = \phi(s_2) = Z_k(s_2) Z_k(s_2+1)^{-1} \\
& \mathfrak{R}_{k_1} = \operatorname*{Res}_{s=1} Z_{k_1}(s),\ \mathfrak{R} = \operatorname*{Res}_{s=1} Z_k(s) \\
& \varrho_{k_1} = \operatorname*{Res}_{s=1} \phi_{k_1}(s) = \frac{\mathfrak{R}_{k_1}}{Z_{k_1}(2)},\ \varrho = \operatorname*{Res}_{s=1} \phi(s) = \frac{\mathfrak{R}}{Z_k(2)} \\
& \mathfrak{V}_{k_1 2} = \varrho_{k_1}^{-1},\ \mathfrak{V}_2 = \varrho^{-1},
\end{aligned}
$$

where $Z_k(s)$ is defined in [3, p. 129] or [5, p. xii]. It is well known that $\mathfrak{V}_{k_1 2}$ (resp. $\mathfrak{V}_2$) is the volume of $G_{1\mathbb{A}}^0/G_{1k}$ (resp. $G_{2\mathbb{A}}^0/G_{2k}$).

Let $C_G = \varrho_{k_1}\varrho$. If $f(w), g(w)$ are meromorphic functions of $w$, we use the notation $f \sim g$ if $f(w) - g(w)$ can be continued meromorphically to a right half plane $\{w \in \mathbb{C} \mid \mathrm{Re}(w) > c\}$ where $c < L(\rho)$ and is holomorphic at $w = L(\rho)$.

Let $J(\Phi, g^0)$ be as in (2.11). We define

$$
(2.17)\qquad I(\Phi, \omega, w) = \int_{G_{\mathbb{A}}^0/G_k} \omega(g^0) J(\Phi, g^0) \mathscr{E}(g^0, w) dg^0.
$$

Then by Shintani's Lemma (3.4.34) [5] and Lemma (2.9),

$$
(2.18)\qquad I(\Phi, \omega, w) \sim C_G \Lambda(w; \rho) I^0(\Phi, \omega).
$$

Now we consider all the paths. Let $\Xi_{\mathfrak{p}}(\Phi, \omega, w)$, etc. be as in Chapter 3 [5]. There are four $\beta$-sequences $\mathfrak{d}_1 = (\beta_1)$, $\mathfrak{d}_2 = (\beta_2)$, $\mathfrak{d}_3 = (\beta_3)$, $\mathfrak{d}_4 = (\beta_2, \beta_3)$. If two paths $\mathfrak{p}_1 = (\mathfrak{d}, \mathfrak{s}_1)$, $\mathfrak{p}_2 = (\mathfrak{d}, \mathfrak{s}_2)$ satisfy $\mathfrak{s}_1(i) = \mathfrak{s}_2(i)$ for all $i > 1$ and $\mathfrak{s}_1(1) = 0$, $\mathfrak{s}_2(1) = 1$, by (2.12),

$$
\Xi_{\mathfrak{p}_2}(\Phi, \omega, w) = \Xi_{\mathfrak{p}_1}(\widehat{\Phi}, \omega^{-1}, w),
$$

etc. So we only have to consider paths $\mathfrak{p} = (\mathfrak{d}, \mathfrak{s})$ such that $\mathfrak{s}(1) = 0$.

Let $\mathfrak{p}_i = (\mathfrak{d}_i, \mathfrak{s}_i)$ be a path such that $\mathfrak{s}_i(1) = 0$ for $i = 1, 2, 3$ (this determines the path). Let $\mathfrak{p}_{41} = (\mathfrak{d}_4, \mathfrak{s}_{41})$ be the path such that $\mathfrak{s}_{41}(1) = 0$, $\mathfrak{s}_{41}(2) = 0$, and $\mathfrak{p}_{42} = (\mathfrak{d}_4, \mathfrak{s}_{42})$ the path such that $\mathfrak{s}_{41}(1) = 0$, $\mathfrak{s}_{41}(2) = 1$. We define

$$\Xi_{\mathrm{st},i}(\Phi, \omega, w) = \int_{G^0_{\mathbb{A}}/G_k} \Theta_{V^{\mathrm{ss}}_{\mathrm{st},i\,k}}(\Phi, g^0) \mathscr{E}(g^0, w) dg^0 \text{ for } i = 1, 2.$$

This is well defined for $\mathrm{Re}(w) \gg 0$ by (3.4.34) [5].

It is not necessary to know $\beta_1, \beta_3$ to describe $\Xi_{\mathfrak{p}_i}(\Phi, \omega, w)$ for $i = 1, 3$. It is easy to see that $\beta_2$ does not depend on the choice of the inner product on $\mathfrak{t}^*$, and $d(\lambda_1, \lambda_2)^{\beta_2} = \lambda_2$. In order to describe distributions associated with paths $\mathfrak{p}_2, \mathfrak{p}_{41}, \mathfrak{p}_{42}$, we have to know $e_{\mathfrak{p}_2 1}(d(1, \lambda_2))$, $e_{\mathfrak{p}_{4i} 1}(d(\lambda_1, \lambda_2))$, $A_{\mathfrak{p}_2 0}, A_{\mathfrak{p}_{4i} 2}$ for $i = 1, 2$. For the definition of $e_{\mathfrak{p}_2}(d(1, \lambda_2))$, etc., see [5, p. 72]. It is easy to see that

$$
\begin{aligned}
(2.19) \qquad & e_{\mathfrak{p}_2 1}(d(1, \lambda_2)) = e_{\mathfrak{p}_{41} 1}(d(\lambda_1, \lambda_2)) = \lambda_2,\ e_{\mathfrak{p}_{42} 1}(d(\lambda_1, \lambda_2)) = \lambda_2^{-1}, \\
& A_{\mathfrak{p}_2 0} = \{d(1, \lambda_2)) \mid \lambda_2 \le 1\},\ A_{\mathfrak{p}_2 1} = \{d(1, \lambda_2)) \mid \lambda_2 \ge 1\}, \\
& A_{\mathfrak{p}_{41} 2} = \{d(\lambda_1, \lambda_2)) \mid \lambda_2 \le 1\},\ A_{\mathfrak{p}_{42} 2} = \{d(\lambda_1, \lambda_2)) \mid \lambda_2 \ge 1\}.
\end{aligned}
$$

It is easy to see that $A'_{\mathfrak{p}_2 0}, A'_{\mathfrak{p}_2 1}$ and $A'_{\mathfrak{p}_{4i} 2}$ for $i = 1, 2$ in [5, p. 73] are the same as $A_{\mathfrak{p}_2 0}, A_{\mathfrak{p}_2 1}$, and $A_{\mathfrak{p}_{4i} 2}$ for $i = 1, 2$. For $\mathfrak{d}_2$,

$$M^1_{\mathfrak{d}_2 \mathbb{A}} = M^0_{\mathfrak{d}_2 \mathbb{A}} = \{(g_1, a(t_{21}, t_{22})) \mid g_1 \in G^0_{1\mathbb{A}},\ t_{11}, t_{12} \in \mathbb{A}^1\}.$$

So if we put $t_2 = a(\underline{\lambda}_2^{-1} t_{21}, \underline{\lambda}_2 t_{22})$,

$$(2.20) \qquad A_{\mathfrak{p}_2 0} M^0_{\mathfrak{d}_2 \mathbb{A}} = \{(g_1, t_2) \in G^0_{1\mathbb{A}} \times T^0_{2\mathbb{A}} \mid \lambda_2 \le 1\}.$$

For $\mathfrak{d}_4$, $M^1_{\mathfrak{d}_4 \mathbb{A}} = M^0_{\mathfrak{d}_4 \mathbb{A}} = \widehat{T}^0_{\mathbb{A}}$. So

$$(2.21) \qquad A_{\mathfrak{p}_{41} 2} M^0_{\mathfrak{d}_4 \mathbb{A}} = \{t^0 \in T^0 \mid \lambda_2 \le 1\},\ A_{\mathfrak{p}_{42} 2} M^0_{\mathfrak{d}_4 \mathbb{A}} = \{t^0 \in T^0 \mid \lambda_2 \ge 1\}.$$

We consider $\Xi_{\mathrm{st},1}(\Phi, \omega, w)$ in §4, $\Xi_{\mathrm{st},2}(\Phi, \omega, w)$ in §5, $\Xi_{\mathfrak{p}_i}(\Phi, \omega, w)$ for $i = 1, 3$ in §6, and $\Xi_{\mathfrak{p}_2}(\Phi, \omega, w)$, $\Xi_{\mathfrak{p}_{41}}(\Phi, \omega, w), \Xi_{\mathfrak{p}_{42}}(\Phi, \omega, w)$ in §7.

By (3.4.34), (3.5.5), (3.5.7) [5] and the consideration in [5, p. 90],

$$
\begin{aligned}
(2.22) \qquad I(\Phi, \omega, w) = {} & \delta_{\#}(\omega) \mathfrak{V}_{2,k_1} \mathfrak{V}_2 (\widehat{\Phi}(0) - \Phi(0)) C_G \Lambda(w; \rho) \\
& + \Xi_{\mathfrak{p}_1}(\widehat{\Phi}, \omega^{-1}, w) + \widetilde{\Xi}_{\mathfrak{p}_1}(\widehat{\Phi}, \omega^{-1}, w) \\
& - (\Xi_{\mathfrak{p}_1}(\Phi, \omega, w) + \widetilde{\Xi}_{\mathfrak{p}_1}(\Phi, \omega, w)) \\
& + \sum_{i=2,3} (\Xi_{\mathfrak{p}_i}(\widehat{\Phi}, \omega^{-1}, w) - \Xi_{\mathfrak{p}_i}(\Phi, \omega, w)) \\
& + \sum_{i=1,2} (\Xi_{\mathrm{st},i}(\widehat{\Phi}, \omega^{-1}, w) - \Xi_{\mathrm{st},i}(\Phi, \omega, w)), \\
\Xi_{\mathfrak{p}_2}(\Phi, \omega, w) = {} & \Xi_{\mathfrak{p}_2 +}(\Phi, \omega, w) + \widehat{\Xi}_{\mathfrak{p}_2 +}(\Phi, \omega, w) \\
& + \Xi_{\mathfrak{p}_2 \#}(\Phi, \omega, w) - \widehat{\Xi}_{\mathfrak{p}_2 \#}(\Phi, \omega, w) \\
& + \Xi_{\mathfrak{p}_{42}}(\Phi, \omega, w) - \Xi_{\mathfrak{p}_{41}}(\Phi, \omega, w),
\end{aligned}
$$

and all the above distributions are well defined.

For the rest of this section, we recall some notations and results. Let

$$(2.23)\qquad \begin{aligned} \mathscr{E}_N(g^0,w) &= \int_{N_{\mathbb{A}}/N_k} \mathscr{E}(g^0 n(u),w)\,du, \\ \mathscr{E}_{N_2}(g^0,w) &= \int_{N_{2\mathbb{A}}/N_{2k}} \mathscr{E}(g^0 n(u_2),w)\,du_2 \end{aligned}$$

be the constant terms of $\mathscr{E}(g^0,w)$ with respect to $N, N_2$ respectively. We define $\widetilde{\mathscr{E}}(g^0,w) = \mathscr{E}(g^0,w) - \mathscr{E}_N(g^0,w)$. For $\alpha = (\alpha_1,\alpha_2) \in k_1 \times k$ and $u = (u_1,u_2) \in \mathbb{A}_{k_1} \times \mathbb{A}$, we define

$$(2.24).\qquad < \alpha n(u) > = < \mathrm{Tr}_{k_1/k}(\alpha_1 u_1) >< \alpha_2 u_2 >$$

Let

$$(2.25)\qquad \mathscr{E}_\alpha(g^0,w) = \int_{N_{\mathbb{A}}/N_k} \mathscr{E}(g^0 n(u),w) < \alpha n(u) > du.$$

Of course, if $\alpha = (0,0)$, $\mathscr{E}_\alpha(g^0,w) = \mathscr{E}_N(g^0,w)$.

Any element of the Weyl group of $G$ is of the form $\tau = (\tau_1,\tau_2)$ where $\tau_1,\tau_2$ are either 1 or $(1,2)$ which corresponds to $\nu$. We define

$$(2.26)\qquad \begin{aligned} & s_\tau = (s_{\tau 1}, s_{\tau 2}) = (z_{1\tau_1(2)} - z_{1\tau_1(1)}, z_{2\tau_2(2)} - z_{2\tau_2(1)}), \\ & M_{\tau_1}(s_{\tau_1}) = \begin{cases} 1 & \tau_1 = 1, \\ \phi_{k_1}(s_{\tau_1}) & \tau_1 = (1,2), \end{cases} \\ & M_{\tau_2}(s_{\tau_2}) = \begin{cases} 1 & \tau_2 = 1, \\ \phi(s_{\tau_2}) & \tau_2 = (1,2), \end{cases} \\ & M_\tau(s_\tau) = M_{\tau_1}(s_{\tau_1}) M_{\tau_2}(s_{\tau_2}). \end{aligned}$$

For each $\tau$, $z$ is determined by $s_\tau$. So any function $f(z)$ of $z$ can be considered as a function of $s_\tau$. We denote this function by $\widetilde{f}(s_\tau)$. For example,

$$(2.27)\qquad L(z) = \widetilde{L}(s_\tau) = \begin{cases} -(s_{\tau 1} + C s_{\tau 2}) & \tau = (1,1), \\ s_{\tau 1} - C s_{\tau 2} & \tau = (\nu,1), \\ -s_{\tau 1} + C s_{\tau 2} & \tau = (1,\nu), \\ s_{\tau 1} + C s_{\tau 2} & \tau = (\nu,\nu). \end{cases}$$

Regarding $\Lambda(w;z)$ as a function of $s_\tau$, we use the notation $\widetilde{\Lambda}(w;s_\tau)$ for this function. We define $\widetilde{\Lambda}_\tau(w;s_\tau) = M_\tau(s_\tau)\widetilde{\Lambda}(w;s_\tau)$. We choose $r = (r_1,r_2) \in \mathbb{R}^2$ so that $r_1, r_2 > 1$. Then it is well known that

$$(2.28)\qquad \mathscr{E}_N(b^1,w) = \sum_\tau \left(\frac{1}{2\pi\sqrt{-1}}\right)^2 \int_{\mathrm{Re}(s_\tau)=(r_1,r_2)} \lambda_1^{s_{\tau 1}-1} \lambda_2^{s_{\tau 2}-1} \widetilde{\Lambda}_\tau(w;s_\tau)\,ds_\tau.$$

We define

$$(2.29) \qquad G(\lambda_1, w) = \frac{\varrho}{2\pi\sqrt{-1}} \int_{\mathrm{Re}(s_{\tau 1})=r_1>1} \lambda_1^{s_{\tau 1}-1} \phi_{k_1}(s_{\tau 1}) \widetilde{\Lambda}(w; s_{\tau 1}, 1) ds_\tau,$$

where $\tau = (\nu, \nu)$ on the right hand side.

**Proposition (2.30)** *Let $M > L(\rho)$. Then there exist constants $\delta, c_1, \cdots, c_k > 0$, $d_0, \cdots, d_k$ such that for any $l_1, l_2 \gg 0$,*

$$(1) \qquad \widetilde{\mathscr{E}}(b^1, w),\ \mathscr{E}_\alpha(b^1, w) \ll \sup(\lambda_1^{l_1}\lambda_2^{d_0}, \lambda_1\lambda_2^{l_2}, \lambda_1^{l_1}\lambda_2^{l_2})$$

*for $\alpha \neq (0,0)$ and $L(\rho) - \delta \leq \mathrm{Re}(w) \leq M$, and*

$$(2) \qquad \mathscr{E}_N(b^1, w) - G(\lambda_1, w) \ll \sum_{i=1}^{k} \lambda_1^{c_i}\lambda_2^{d_i}.$$

*Proof.* Consider $b^1 = (b_1, b_2)$ in (2.5). Let $E_{N_i}(b_i, z_i)$, $\widetilde{E}_{B_i}(b_i, z_i)$ be the constant term and the non-constant term of $E_{B_i}(b_i, z_i)$ for $i = 1, 2$. Then $\widetilde{\mathscr{E}}(b^1, w)$ is a sum of contour integrals of

$$\widetilde{E}_{B_1}(b_1, z_1)E_{N_2}(b_2, z_2),\ E_{N_1}(b_1, z_1)\widetilde{E}_{B_2}(b_2, z_2),\ \widetilde{E}_{B_1}(b_1, z_1)\widetilde{E}_{B_2}(b_2, z_2).$$

On any vertical strip in $\mathrm{Re}(z_{i1} - z_{i2}) > 0$ and for $l_i \gg 0$, $\widetilde{E}_{B_i}(b_i, z_i) \ll \lambda_i^{l_i}$ for $i = 1, 2$. Moreover, this estimate applies to each term of the Fourier expansion.

Note that

$$E_{N_i}(b_i, z_i) = \sum_{\tau_i} M_{\tau_i}(s_{\tau_i i})\lambda_i^{\tau_i - 1},$$

where $s_{\tau_i i}$ is defined as in (2.26) (we are only considering each factor of the Weyl group). We divide the contour integral of $\widetilde{E}_{B_1}(b_1, z_1)E_{N_2}(b_2, z_2)$ according to $\tau_2$ and choose the contour so that $\mathrm{Re}(z_{11} - z_{12}) = \frac{1}{2}$ and $\mathrm{Re}(s_{\tau 2})$ close to 1. This corresponds to the first term in (1). For $E_{N_1}(b_1, z_1)\widetilde{E}_{B_2}(b_2, z_2)$, we choose the contour so that $\mathrm{Re}(s_{\tau_1 1}) = 2$, $\mathrm{Re}(z_{21} - z_{22}) = \frac{1}{2}$. Then $\pm 2 + \frac{C}{2} < 1 + C = L(\rho)$, because $C > 4$ by assumption. This corresponds to the second term in (1). For $\widetilde{E}_{B_1}(b_1, z_1)\widetilde{E}_{B_2}(b_2, z_2)$, we choose $\mathrm{Re}(z_{i1} - z_{i2}) = \frac{1}{2}$ for $i = 1, 2$. This corresponds to the third term in (1). This proves (1).

Consider (2). Let $\tau$ be a Weyl group element. If $\tau = (1, 1)$ or $(\nu, 1)$, we choose $\mathrm{Re}(s_\tau) = (2, \frac{1}{2})$. If $\tau = (1, \nu)$, we choose $\mathrm{Re}(s_\tau) = (2, r_2)$ where $r_2$ is close to 1. Then $\widetilde{L}(2, r_2) = -2 + Cr_2 < 1 + C$. So the corresponding terms in $\mathscr{E}_N(b^1, w)$ satisfy the estimate of (2). Let $\tau = (\nu, \nu)$. We move the contour to $\mathrm{Re}(s_\tau) = (1 + \delta_1, 1 - \delta_2)$ where $\delta_1^{-1}\delta_2 \gg 0$. Then $1 + \delta_1 + C(1 - \delta_2) < 1 + C$. When the contour crosses the line $s_{\tau 2} = 1$, we get the residue $G(\lambda_1, w)$. The contour integral over $\mathrm{Re}(s_\tau) = (1 + \delta_1, 1 - \delta_2)$ satisfies the estimate of (2). This proves the proposition.

Q.E.D.

In particular, if we are not interested in the growth, $\widetilde{\mathscr{E}}(b^1, w)$, $\mathscr{E}_\alpha(b^1, w)$ are bounded by a slowly increasing function (but across the point $L(\rho)$) and this bound is uniform with respect to $\alpha \neq (0,0)$.

## §3 The zeta function for the space of binary Hermitian forms

Since the space $W$ of (single) binary Hermitian forms appears in the induction, we have to know the principal part formula of the zeta function for this case. This case is quite easy and is well known. The purpose of this section is to review the principal part formula for this case. This is a case where the theta series for the entire lattice $V_k$ is integrable on $G^0_{\mathbb{A}}$.

Before discussing the space of binary Hermitian forms, we make a few definitions for later purposes.

If $\Psi \in \mathscr{S}(\mathbb{A})$, $\omega$ is a character of $\mathbb{A}^\times/k^\times$, $t \in \mathbb{A}^\times$, and $s \in \mathbb{C}$, we define

$$
\begin{aligned}
\Theta_1(\Psi, t) &= \sum_{x \in k^\times} \Psi(tx), \\
\Sigma_1(\Psi, \omega, s) &= \int_{\mathbb{A}^\times/k^\times} \omega(t)|t|^s d^\times t.
\end{aligned} \tag{3.1}
$$

Properties of $\Sigma_1(\Psi, \omega, s)$ are well known.

If $\Psi \in \mathscr{S}(\mathbb{A} \times \mathbb{A}_{k_1})$, $\omega = (\omega_1, \omega_2)$ is a character of $\mathbb{A}^\times/k^\times \times \mathbb{A}^\times_{k_1}/k_1^\times$, $t = (t_1, t_2) \in \mathbb{A}^\times \times \mathbb{A}^\times_{k_1}$, and $s = (s_1, s_2) \in \mathbb{C}^2$, we define

$$
\begin{aligned}
\Theta_{1,1}(\Psi, t) &= \sum_{x \in k^\times \times k_1^\times} \Psi(tx), \\
\Sigma_{1,1}(\Psi, \omega, s) &= \int_{\mathbb{A}^\times/k^\times \times \mathbb{A}^\times_{k_1}/k_1^\times} \omega_1(t_1)\omega_2(t_2)|t_1|^{s_1}|t_2|_{k_1}^{s_2} d^\times t_2 d^\times t_2.
\end{aligned} \tag{3.2}
$$

The function $\Sigma_{1,1}(\Psi, \omega, s)$ is more or less the product $Z_k(s_1)Z_{k_1}(s_2)$ and it is easy to see that it can be continued meromorphically everywhere having poles at $s_1 = 0, 1,\ s_2 = 0, 1$. If $s_0 \in \mathbb{C}$ or $s_0 = (s_{01}, s_{02}) \in \mathbb{C}^2$, let

$$
\begin{aligned}
\Sigma_1(\Psi, \omega, s) &= \sum_{i=-1}^{\infty} \Sigma_{1,(i)}(\Psi, \omega, s_0)(s - s_0)^i, \\
\Sigma_{1,1}(\Psi, \omega, s) &= \sum_{i,j=-1}^{\infty} \Sigma_{1,1,(i,j)}(\Psi, \omega, s_0)(s_1 - s_{01})^i(s_2 - s_{02})^j
\end{aligned} \tag{3.3}
$$

be the Laurent expansions around $s = s_0$. If $\omega$ is trivial, we drop $\omega$ and use notations $\Sigma_1(\Psi, s), \Sigma_{1,(i)}(\Psi, s_0)$, etc.

Let $i$ be a positive integer. For $\mu = (\mu_1, \cdots, \mu_i) \in \mathbb{R}^i_+$ and $M > 0$, we define

$$
\mathrm{rd}_{i,M}(\mu) = \inf(\mu_1^{\pm 1} \cdots \mu_i^{\pm 1}), \tag{3.4}
$$

where we consider all the possible $\pm$.

Let $W$ be as in §1 and $G = \mathrm{GL}(1)_k \times \mathrm{GL}(2)_{k_1}$ considered as a group over $k$. The action of $\mathrm{GL}(2)_{k_1}$ is the same as in §1 and $\alpha \in \mathrm{GL}(1)_k$ acts by multiplication by $\alpha$. Let $G_1 = \mathrm{GL}(1)_k,\ G_2 = \mathrm{GL}(2)_{k_1}$. We define $G^0_{\mathbb{A}},\ G^0_{1\mathbb{A}},\ \widetilde{G}_{\mathbb{A}}$, etc. similarly as in §2. If $\omega = (\omega_1, \omega_2)$ is a character of $\mathbb{A}_{k_1}^\times / k_1^\times \times \mathbb{A}^\times / k^\times$, we define $\omega(t_1, g_2) = \omega_1(t_1)\omega_2(\det g_2)$ for $t_1 \in \mathbb{A}^1$, $g_2 \in G_{2\mathbb{A}}^0$.

We write elements of $W$ as $x = \begin{pmatrix} x_0 & x_1 \\ x_1^\sigma & x_2 \end{pmatrix}$. It is easy to see that $W^{\mathrm{ss}}$ is the set of non-singular matrices. Let $S_{W,1}$ be the set of rank one matrices in $W$. Let

$$Z_{W,1} = \left\{ x = \begin{pmatrix} 0 & 0 \\ 0 & x_2 \end{pmatrix} \right\},\ Z_{W,1}^{\mathrm{ss}} = \{x \in Z_{W,1} \mid x_2 \neq 0\}.$$

Let $B = \mathrm{GL}(1)_k \times B_2 \subset G$ be the Borel subgroup as in §2. Then it is easy to see $S_{W,1k} \cong G_k \times_{B_k} Z_{W,1k}^{\mathrm{ss}}$.

Let $\Phi \in \mathscr{S}(W_{\mathbb{A}})$. For any $G_k$–invariant subset $S \subset W_k$, we define $\Theta_S(\Phi, \widetilde{g})$ in the usual manner for $\widetilde{g} \in \widetilde{G}_{\mathbb{A}}$. For $g_2 \in G_{2\mathbb{A}}^0$, let $g_2 = ka(\mu^{-1}, \mu)b$ be the Iwasawa decomposition where $b \in \widehat{T}_{2\mathbb{A}}^0 N_{2\mathbb{A}}$.

If $x \in W_k^{\mathrm{ss}}$, then $x_0 \neq 0$, or $x_0 = 0$, $x_1 \neq 0$. So by Lemma (1.2.6) [5], there exists a slowly increasing function $h(\lambda, \mu)$ such that for any $N_1, N_2 \geq 1$,

$$\Theta_{W_k^{\mathrm{ss}}}(\Phi, \underline{\lambda} g^0) \ll \sup(h(\lambda, \mu)(\lambda\mu^{-1})^{-N_1}, \lambda^{-2N_1} \sup(1, \mu^{-1})) \tag{3.5}$$

for $g^0$ in a Siegel set and $\lambda \in \mathbb{R}_+$. Also by Lemma (1.2.6) [5],

$$\Theta_{W_k}(\Phi, g^0) \ll \sup(1, \mu) \sup(1, \mu^{-1}) \ll \mu^{-1}$$

for $g^0$ in a Siegel set. This implies that $\Theta_{W_k}(\Phi, g^0)$ is integrable on $G_{\mathbb{A}}^0 / G_k$.

We define

$$Z_W(\Phi, \omega, s) = \int_{\widetilde{G}/G_k} \omega(\widetilde{g}) \lambda^s \Theta_{W_k^{\mathrm{ss}}}(\Phi, \widetilde{g}) d\widetilde{g},$$

$$Z_{W+}(\Phi, \omega, s) = \int_{\substack{\widetilde{G}/G_k \\ \lambda \geq 1}} \omega(\widetilde{g}) \lambda^s \Theta_{W_k^{\mathrm{ss}}}(\Phi, \widetilde{g}) d\widetilde{g}.$$

Then $Z_W(\Phi, \omega, s)$ is defined on a certain right half plane and $Z_{W+}(\Phi, \omega, s)$ is an entire function.

For $x, y \in W$, we define $[x, y]' = \mathrm{Tr}_{k_1/k}(\mathrm{tr}(xy))$, $[x, y] = [x, (1, \nu)y]'$. We use this non-degenerate bilinear form to define a Fourier transform $\widehat{\Phi}$ of $\Phi$. We define $\widetilde{g}^\iota$ similarly as in §2.

Let $R_{W,1}\Phi$ be the restriction of $\Phi$ to $Z_{W,1\mathbb{A}}$. We put $\delta_{W,1}(\omega) = \delta(\omega_1 \circ \mathrm{N}_{k_1/k})\delta(\omega_2)$ where $\omega_1 \circ \mathrm{N}_{k_1/k}$ is the composition of $\omega_1, \mathrm{N}_{k_1/k}$. We define $M_\omega \Phi$ as in (3.1.11) [5] and assume that $\Phi = M_\omega \Phi$. Since there is no difficulty regarding the convergence of the full theta series $\Theta_{W_k}(\Phi, g^0)$, the following formula can easily be proved and the details are left to the reader.

$$\begin{aligned} Z_W(\Phi, \omega, s) &= Z_{W+}(\Phi, \omega, s) + Z_{W+}(\widehat{\Phi}, \omega^{-1}, 4 - s) \\ &+ \delta_\#(\omega) \mathfrak{V}_{2,k_1} \left( \frac{\widehat{\Phi}(0)}{s-4} - \frac{\Phi(0)}{s} \right) \\ &+ \delta_{W,1}(\omega) \left( \frac{\Sigma_1(R_{W,1}\widehat{\Phi}, \omega_1^{-1}, 2)}{s-2} - \frac{\Sigma_1(R_{W,1}\Phi, \omega_1, 2)}{s-2} \right). \end{aligned} \tag{3.6}$$

If $s_0 \in \mathbb{C}$, let

$$Z_W(\Phi,\omega,s) = \sum_{i=-1}^{\infty} Z_{W,(i)}(\Phi,\omega,s_0)(s-s_0)^i \tag{3.7}$$

be the Laurent expansion around $s = s_0$.

## §4 The adjusting term I

Before we start the analysis of $\Xi_{\mathrm{st},1}(\Phi,\omega,w)$, etc., we introduce some notations. For $\Phi \in \mathscr{S}(V_{\mathbb{A}})$, let $R_i\Phi$ be the restriction of $\Phi$ to $Z_{i\mathbb{A}}$ for $i = 2,3$. Let $Z_4 = Z_3$, $Z_4^{\mathrm{ss}} = Z_3^{\mathrm{ss}}$. For any subspace $Z_4 \subset U \subset V$ and $\Psi \in \mathscr{S}(U_{\mathbb{A}})$, let $R_4\Psi$ be the restriction of $\Psi$ to $Z_{4\mathbb{A}}$. Let $\widetilde{R}_1\Phi$ be the restriction of $\Phi$ to $Y_{1\mathbb{A}}$. Let $\widetilde{R}_{0,1}\Phi$ be the restriction of $\Phi$ to $Y_{0,1\mathbb{A}}$. Let $R_{0,2}\Phi$ be the restriction of $\Phi$ to $Z_{0,2\mathbb{A}}$. We define $R_1\Phi \in \mathscr{S}(Z_{1\mathbb{A}})$, $R_{0,1}\Phi \in \mathscr{S}(Z_{0,1\mathbb{A}})$ by

$$\begin{aligned} R_1\Phi(x_{12},x_{21}) &= \int_{\mathbb{A}} \widetilde{R}_1\Phi(x_{12},x_{21},x_{22})dx_{22}, \\ R_{0,1}\Phi(x_{12},x_{20},x_{21}) &= \int_{\mathbb{A}} \widetilde{R}_{0,1}\Phi(x_{12},x_{20},x_{21},x_{22})dx_{22}. \end{aligned} \tag{4.1}$$

Now we consider $\Xi_{\mathrm{st},1}(\Phi,\omega,w)$.

Let $\Phi$ and $\omega$ be as before. Let $\mu = (\mu_1,\mu_2) \in \mathbb{R}_+^2$, $q = (q_1,q_2) \in (\mathbb{A}^1)^2$, $u_0 \in \mathbb{A}_{k_1}$, and $s_1, s_2 \in \mathbb{C}$. Let $\alpha_{k_1}(u_0)$ be the function of $u_0 \in \mathbb{A}_{k_1}$ defined in Definition (2.2.2) [5] for the field $k_1$. For $\omega = (\omega_1,\omega_2)$ as before, we define $\delta_{\mathrm{st},1}(\omega) = \delta(\omega_1(\omega_2 \circ \mathrm{N}_{k_1/k})^{-1})$. We define

$$\begin{aligned} f_1(\Phi,\omega,\mu,q,u_0,s_1,s_2) &= \delta_{\mathrm{st},1}(\omega)\omega_2(q_1q_2)\mu_2^{s_1}\alpha_{k_1}(u_0)^{s_2} \\ &\quad \times R_{0,1}\Phi(\underline{\mu}_1\underline{\mu}_2^{-1}q_1, \underline{\mu}_1\underline{\mu}_2 q_2, \underline{\mu}_1\underline{\mu}_2 q_2 u_0). \end{aligned} \tag{4.2}$$

We define $d^\times\mu = d^\times\mu_1 d^\times\mu_2$, $d^\times q = d^\times q_1 d^\times q_2$. Let $s \in \mathbb{C}$ be another variable.

**Definition (4.3)**

(1) $T_{\mathrm{st},1}(\Phi,\omega,s,s_1,s_2) = \displaystyle\int_{\mathbb{R}_+^2 \times (\mathbb{A}^1)^2 \times \mathbb{A}_{k_1}} \mu_1^s f_1(\Phi,\omega,\mu,q,u_0,s_1,s_2)d^\times\mu d^\times q du_0,$

(2) $T_{\mathrm{st},1+}(\Phi,\omega,s,s_1,s_2) = \displaystyle\int_{\substack{\mathbb{R}_+^2 \times (\mathbb{A}^1)^2 \times \mathbb{A}_{k_1} \\ \mu_1 \geq 1}} \mu_1^s f_1(\Phi,\omega,\mu,q,u_0,s_1,s_2)d^\times\mu d^\times q du_0,$

(3) $T_{\mathrm{st},1}^1(\Phi,\omega,s_1,s_2) = \displaystyle\int_{\mathbb{R}_+ \times (\mathbb{A}^1)^2 \times \mathbb{A}_{k_1}} f_1(\Phi,\omega,1,\mu_2,q,u_0,s_1,s_2)d^\times\mu d^\times q du_0.$

**Proposition (4.4)** *The integral* (4.3)(1) *converges absolutely and locally uniformly for* $\mathrm{Re}(s) - \mathrm{Re}(s_1) > 2$, $\mathrm{Re}(s) + \mathrm{Re}(s_1) > 6$, $\mathrm{Re}(s) + \mathrm{Re}(s_1) + 4\mathrm{Re}(s_2) > 6$, *and the integrals* (4.3)(2), (4.3)(3) *are entire functions.*

*Proof.* Let $\Psi = R_{0,1}\Phi$. We define

$$(4.5)\qquad \widetilde{\Theta}_{Z_{0,1}}(\Psi,\mu,q,u_0) = \sum_{x,y\in k^\times} \Psi(\underline{\mu}_1\underline{\mu}_2^{-1}q_1x, \underline{\mu}_1\underline{\mu}_2q_2y, \underline{\mu}_1\underline{\mu}_2q_2yu_0).$$

Then it is easy to see that

$$(4.6) T_{\mathrm{st},1}(\Psi,\omega,s,s_1,s_2) = \int_{\mathbb{R}_+^2\times(\mathbb{A}^1/k^\times)^2\times\mathbb{A}_{k_1}} \omega_2(q_1q_2)\mu_1^s\mu_2^{s_1}\alpha_{k_1}(u_0)^{s_2}\widetilde{\Theta}_{Z_{0,1}}(\Psi,\mu,q,u_0)d^\times\mu d^\times q du_0,$$

etc.

Let $\sigma = \mathrm{Re}(s)$, $\sigma_1 = \mathrm{Re}(s_1)$, $\sigma_2 = \mathrm{Re}(s_2)$. By Lemma (1.2.3) [5], there exists $0 \le \Psi_1 \in \mathscr{S}(Z_{0,1\mathbb{A}})$ such that (4.6) is bounded by a constant multiple of

$$(4.7)\qquad \int_{\mathbb{R}_+^2\times\mathbb{A}_{k_1}} \mu_1^\sigma\mu_2^{\sigma_1}\alpha_{k_1}(u_0)^{\sigma_2}\widetilde{\Theta}_{Z_{0,1}}(\Psi_1,\mu,1,u_0)d^\times\mu du_0.$$

We estimate $\alpha_{k_1}(u_0)^{\sigma_2}\widetilde{\Theta}_{Z_{0,1}}(\Psi_1,\mu,1,u_0)$, which is

$$(4.8)\qquad \alpha_{k_1}(u_0)^{\sigma_2}\sum_{x,y\in k^\times}\Psi_1(\underline{\mu}_1\underline{\mu}_2^{-1}x, \underline{\mu}_1\underline{\mu}_2y, \underline{\mu}_1\underline{\mu}_2yu_0).$$

We put $u_3 = \underline{\mu}_1\underline{\mu}_2yu_0$. Then there exists a lattice $L \subset k$ and a compact set $F \subset \mathbb{A}_{k_1 f}$ (the finite part of $\mathbb{A}_{k_1}$) such that

$$\Psi_1(\underline{\mu}_1\underline{\mu}_2^{-1}x, \underline{\mu}_1\underline{\mu}_2y, u_3) = 0$$

unless $x, y \in L$ and the finite part of $u_3$ belongs to $F$.

Consider the function $\alpha_{k_1}(u_0)$. Let $\mathfrak{M}_{k_1\infty}$ be the set of infinite places of $k_1$, and $L \subset L_1 \subset k_1$ a lattice. We proved in [5, p. 117] that if $\sigma_2 \le 0$ and $\lambda \in \mathbb{R}_+$,

$$(4.9)\qquad \alpha_{k_1}(\underline{\lambda}_{k_1}^{-1}y^{-1}u_3)^{\sigma_2} \le \lambda^{\sigma_2}\prod_{v\in\mathfrak{M}_{k_1\infty}}\alpha_v(\underline{\lambda}_{k_1}y)^{\sigma_2}\alpha_v(u_3)^{\sigma_2}$$

If $\sigma_2 \ge 0$, $\alpha_{k_1}(u_3)^{\sigma_2} \le 1$. Note that $\prod_{v\in\mathfrak{M}_{k_1\infty}}\alpha_v(\underline{\lambda}_{k_1}y)^{\sigma_2}\alpha_v(u_3)^{\sigma_2}$ is a polynomial growth function of the infinite part of $\underline{\lambda}y$, $u_3$ as long as $\sigma$ ranges over a compact set. We apply (4.9) to $\lambda = (\mu_1\mu_2)^2$. Since $\Psi_1$ is a Schwartz–Bruhat function, there exists $0 \le \Psi_2 \in \mathscr{S}(Z_{0,1\mathbb{A}})$ such that (4.8) is bounded by a constant multiple of

$$(4.10)\qquad \sup(1,(\mu_1\mu_2)^{2\sigma_2})\sum_{x,y\in k^\times}\Psi_2(\underline{\mu}_1\underline{\mu}_2^{-1}x, \underline{\mu}_1\underline{\mu}_2y, u_3).$$

Let

$$\Psi_3(x,y) = \int_{\mathbb{A}_{k_1}}\Psi_2(x,y,z)dz.$$

We choose $|\Psi_3| \le \Psi_4 \in \mathscr{S}(\mathbb{A}^2)$. Then (4.7) is bounded by a constant multiple of

$$\int_{\mathbb{R}_+^2} \sup(1, (\mu_1\mu_2)^{2\sigma_2})(\mu_1\mu_2)^{-2}\mu_1^{\sigma}\mu_2^{\sigma_1} \sum_{x,y\in k^\times} \Psi_3(\underline{\mu}_1\underline{\mu}_2^{-1}x, \underline{\mu}_1\underline{\mu}_2 y)d^\times\mu. \tag{4.11}$$

Since $\mu_1^{\sigma}\mu_2^{\sigma_1} = (\mu_1\mu_2^{-1})^{\frac{\sigma-\sigma_1}{2}}(\mu_1\mu_2)^{\frac{\sigma+\sigma_1}{2}}$, (4.11) converges absolutely if

$$\frac{\sigma-\sigma_1}{2} > 1,\ \frac{\sigma+\sigma_1}{2} - 2 > 1,\ \frac{\sigma+\sigma_1}{2} + 2\sigma_2 - 2 > 1.$$

If $\mu_1 \ge 1$ or $\mu_1 = 1$, the integral analogous to (4.11) always converges. This proves Proposition (4.4).

Q.E.D.

Let

$$b_{\mathrm{st},1} = (n(u_0), 1)t^0(1, n(u_2)),$$

where $u_0 \in \mathbb{A}_{k_1}$, $u_2 \in \mathbb{A}$, $t^0 \in T^0_{\mathbb{A}}$. We define

$$X_1 = \{b_{\mathrm{st},1} \mid \lambda_1 \le \alpha_{k_1}(u_0)^{-\frac{1}{2}}\}. \tag{4.12}$$

Then as in [5, p. 112], $KX_1/T_k$ is the fundamental domain for $G^0_{\mathbb{A}}/H_{1k}$. Note that the right invariant measure on $X_1$ is $db_{\mathrm{st},1} = \lambda_2^2 d^\times t^0 du_0 du_2$.

**Proposition (4.13)**

$$\Xi_{\mathrm{st},1}(\Phi, \omega, w) \sim \int_{X_1/T_k} \Theta_{Z'^{\mathrm{ss}}_{0,1k}}(\Phi, b_{\mathrm{st},1})\mathscr{E}_N(b_{\mathrm{st},1}, w)db_{\mathrm{st},1}.$$

*Proof.* By (1.2),

$$\Xi_{\mathrm{st},1}(\Phi, \omega, w) = \int_{X_1/T_k} \Theta_{Z'^{\mathrm{ss}}_{0,1k}}(\Phi, b_{\mathrm{st},1})\mathscr{E}(b_{\mathrm{st},1}, w)db_{\mathrm{st},1}.$$

By Proposition (2.30), we only have to prove that for a fixed slowly increasing function $h(\lambda_1, \lambda_2)$, the integral

$$\int_{X_1/T_k} \Theta_{Z'^{\mathrm{ss}}_{0,1k}}(\Phi, b_{\mathrm{st},1})h(\lambda_1, \lambda_2)\sup(\lambda_1^{l_1}, \lambda_2^{l_2}, \lambda_1^{l_1}\lambda_2^{l_2})db_{\mathrm{st},1} \tag{4.14}$$

converges absolutely for $l_1, l_2 \gg 0$.

In the above integral, $\Theta_{Z'^{\mathrm{ss}}_{0,1k}}(\Phi, b_{\mathrm{st},1})$ is the only part which depends on $u_0, u_2$. It is easy to see that

$$\int_{\mathbb{A}} \Theta_{Z'^{\mathrm{ss}}_{0,1k}}(\Phi, b_{\mathrm{st},1})du_2 = (\lambda_1\lambda_2)^{-1}\widetilde{\Theta}_{Z_{0,1}}(\Psi, \mu, q, u_0), \tag{4.15}$$

where $\Psi = R_{0,1}\Phi$, $\mu_1 = 1$, $\mu_2 = \lambda_1^{-1}\lambda_2$, $q_1 = \mathrm{N}_{k_1}(t_{12})t_{21}$, and $q_2 = \mathrm{N}_{k_1}(t_{10})t_{22}$. Since $\lambda_2 = \lambda_1\mu_2$,

$$\lambda_2^{l_2} = \lambda_1^{l_2}\mu_2^{l_2},\ \lambda_1^{l_1}\lambda_2^{l_2} = \lambda_1^{l_1+l_2}\mu_2^{l_2}.$$

We choose $|\Psi| \leq \Psi_1 \in \mathscr{S}(Z_{0,1\mathbb{A}})$. Then (4.13) is bounded by a constant multiple of a finite sum of integrals of the form

$$\int_{\substack{\mathbb{R}_+^2 \times (\mathbb{A}^1/k^\times)^2 \times \mathbb{A}_{k_1} \\ \lambda_1 \leq \alpha_{k_1}(u_0)^{-\frac{1}{2}}}} \lambda_1^{c_1} \mu_2^{c_2} \widetilde{\Theta}_{Z_{0,1}}(\Psi_1, 1, \mu_2, q, u_0) d^\times \lambda_1 d^\times \mu_2 d^\times q du_0,$$

where $c_1 > 0$. Then the convergence of this integral follows from (4.4) (for (4.3)(3)).

Q.E.D.

**Proposition (4.16)** *Let $\tau = (\nu, \nu)$. Then,*

$$\begin{aligned} &\Xi_{\mathrm{st},1}(\Phi, \omega, w) \\ &\sim \frac{\varrho}{2\pi\sqrt{-1}} \int_{\mathrm{Re}(s_{\tau 1}) = r_1 > 1} \frac{T^1_{\mathrm{st},1}(R_{0,1}\Phi, \omega, 1, -\frac{s_{\tau 1}-1}{2})}{s_{\tau 1} - 1} \phi_{k_1}(s_{\tau 1}) \widetilde{\Lambda}(w; s_{\tau 1}, 1) ds_{\tau 1}. \end{aligned}$$

*Proof.* If $f(q)$ is a function of $q = (q_1, q_2) \in (\mathbb{A}^1/k^\times)^2$,

$$\begin{aligned} (4.17) \qquad &\int_{\widehat{T}^0_{\mathbb{A}}/T_k} \omega(\widehat{t}^0) f(\mathrm{N}_{k_1/k}(t_{12}) t_{21}, \mathrm{N}_{k_1/k}(t_{11}) t_{22}) d^\times \widehat{t}^0 \\ &= \delta(\omega_1(\omega_2 \circ \mathrm{N}_{k_1/k})^{-1}) \int_{(\mathbb{A}^1/k^\times)^2} \omega_2(q_1 q_2) f(q_1, q_2) d^\times q_1 d^\times q_2. \end{aligned}$$

Let $\Psi = R_{0,1}\Phi$, $\mu_2 = \lambda_1^{-1}\lambda_2$. Then $\lambda_2 = \lambda_1\mu_2$ and $d^\times\lambda_1 d^\times\lambda_2 = d^\times\lambda_1 d^\times\mu_2$. So

$$\begin{aligned} (4.18) \qquad &\mathscr{E}_N(b_{\mathrm{st},1}, w) \\ &= \sum_\tau \left(\frac{1}{2\pi\sqrt{-1}}\right)^2 \int_{\mathrm{Re}(s_\tau) = (r_1, r_2)} \lambda_1^{s_{\tau 1} + s_{\tau 2} - 2} \mu_2^{s_{\tau 2} - 1} \widetilde{\Lambda}_\tau(w; s_\tau) ds_\tau. \end{aligned}$$

Note that $\lambda_2^2(\lambda_1\lambda_2)^{-1} = \lambda_1^{-1}\lambda_2 = \mu_2$. By (4.13), (4.15), (4.17),

$$\begin{aligned} \Xi_{\mathrm{st},1}(\Phi, \omega, w) &\sim \delta_{\mathrm{st},1}(\omega) \int_{\substack{\mathbb{R}_+^2 \times (\mathbb{A}^1/k^\times)^2 \times \mathbb{A}_{k_1} \\ \lambda_1 \leq \alpha_{k_1}(u_0)^{-\frac{1}{2}}}} \omega_2(q_1 q_2) \mu_2 \lambda_1^{s_{\tau 1} + s_{\tau 2} - 2} \mu_2^{s_{\tau 2} - 1} \\ &\qquad \times \widetilde{\Theta}_{Z_{0,1}}(\Psi, \mu, q, u_0) d^\times\lambda_1 d^\times\mu_2 d^\times q du_0 \\ &= \delta_{\mathrm{st},1}(\omega) \int_{\mathbb{R}_+ \times (\mathbb{A}^1/k^\times)^2 \times \mathbb{A}_{k_1}} \omega_2(q_1 q_2) \mu_2^{s_{\tau 2}} \frac{\alpha_{k_1}(u_0)^{-\frac{s_{\tau 1} + s_{\tau 2} - 2}{2}}}{s_{\tau 1} + s_{\tau 2} - 2} \\ &\qquad \times \widetilde{\Theta}_{Z_{0,1}}(\Psi, \mu, q, u_0) d^\times\mu_2 d^\times q du_0 \\ &= \frac{T^1_{\mathrm{st},1}(\Psi, \omega, s_{\tau 2}, -\frac{s_{\tau 1} + s_{\tau 2} - 2}{2})}{s_{\tau 1} + s_{\tau 2} - 2}. \end{aligned}$$

So

$$\Xi_{\mathrm{st},1}(\Phi, \omega, w) \sim \left(\frac{1}{2\pi\sqrt{-1}}\right)^2 \int_{\mathrm{Re}(s_\tau) = (r_1, r_2)} \frac{T^1_{\mathrm{st},1}(\Psi, \omega, s_{\tau 2}, -\frac{s_{\tau 1} + s_{\tau 2} - 2}{2})}{s_{\tau 1} + s_{\tau 2} - 2} ds_\tau.$$

If $\tau \neq (\nu, \nu)$, we choose $r_1, r_2 > 1$ close enough to 1. This implies $\widetilde{L}(r) < L(\rho)$. So we only have to consider $\tau = (\nu, \nu)$. Then we get the proposition by moving the contour from $r_1, r_2 > 1$ to $(2, \frac{1}{2})$ as before.

Q.E.D.

We define the adjusting term corresponding to $V_{\mathrm{st},1k}^{\mathrm{ss}}$ as follows.

**Definition (4.19)**

(1) $$T_{\mathrm{st},1}(\Phi, \omega, s, s_1) = \left.\frac{d}{ds_2}\right|_{s_2=0} T_{\mathrm{st},1}(\Phi, \omega, s, s_1, s_2),$$

(2) $$T_{\mathrm{st},1+}(\Phi, \omega, s, s_1) = \left.\frac{d}{ds_2}\right|_{s_2=0} T_{\mathrm{st},1+}(\Phi, \omega, s, s_1, s_2),$$

(3) $$T_{\mathrm{st},1}^1(\Phi, \omega, s_1) = \left.\frac{d}{ds_2}\right|_{s_2=0} T_{\mathrm{st},1}^1(\Phi, \omega, s_1, s_2).$$

By Proposition (4.4), $T_{\mathrm{st},1}(\Phi, \omega, s, s_1)$ is holomorphic for $\mathrm{Re}(s) + \mathrm{Re}(s_1) > 6$, and $T_{\mathrm{st},1+}(\Phi, \omega, s, s_1)$, $T_{\mathrm{st},1}^1(\Phi, \omega, s_1)$ are entire functions.

## §5 The adjusting term II

We consider $\Xi_{\mathrm{st},2}(\Phi, \omega, w)$ in this section.

Let

$$g_{\mathrm{st},2} = (n(u_0)t_1, g_2),\ g_{\mathrm{st},2}^0 = (n(u_0)a(t_{11}, t_{12}), g_2),$$

where $u_0 \in \mathbb{A}_{k_1}$, $t_1 = a(t_{11}\underline{\lambda}_{1k_1}^{-1}, t_{12}\underline{\lambda}_{1k_1}) \in T_{1\mathbb{A}}^0$, $t_{11}, t_{12} \in \mathbb{A}_{k_1}^1$, and $g_2 \in G_{2\mathbb{A}}^0$. We define

(5.1) $$X_2 = \{g_{\mathrm{st},2} \mid \lambda_1 \leq \alpha_{k_1}(u_0)^{-\frac{1}{2}}\}.$$

Let

$$dg_{\mathrm{st}} = d^\times\lambda_1 d^\times t_{11} d^\times t_{12} du_0 dg_2,\ dg_{\mathrm{st}}^0 = d^\times t_{11} d^\times t_{12} du_0 dg_2.$$

Then $dg_{\mathrm{st}}$ is an invariant measure on $X_2$. As in [5, p. 112], $K_1 X_2/(T_{1k} \times G_{2k})$ is the fundamental domain for $G_{\mathbb{A}}^0/H_{2k}$.

For $x, y \in k_1$, we define

(5.2) $$w_{x,y} = v_1 \begin{pmatrix} 0 & x \\ x^\sigma & 0 \end{pmatrix} + v_2 \begin{pmatrix} 0 & y \\ y^\sigma & 0 \end{pmatrix}.$$

Let $s, s_1 \in \mathbb{C}$. We define

(5.3) $$f_2(\Phi, \lambda, g_{\mathrm{st},2}^0, s, s_1) = \omega(g_{\mathrm{st},2}^0)\lambda^s \alpha_{k_1}(u_0)^{s_1} \Theta_{Z'^{\mathrm{ss}}_{0,2k}}(\Phi, \underline{\lambda} g_{\mathrm{st},2}^0).$$

**Definition (5.4)**

(1) $$T_{\mathrm{st},2}(\Phi, \omega, s, s_1) = \int_{\mathbb{R}_+ \times (\mathbb{A}_{k_1}^1/k_1^\times)^2 \times \mathbb{A}_{k_1} \times G_{2\mathbb{A}}^0/G_{2k}} f_2(\Phi, \lambda, g_{\mathrm{st},2}^0, s, s_1) d^\times\lambda dg_{\mathrm{st},2}^0,$$

(2) $$T_{\mathrm{st},2+}(\Phi, \omega, s, s_1) = \int_{\substack{\mathbb{R}_+ \times (\mathbb{A}_{k_1}^1/k_1^\times)^2 \times \mathbb{A}_{k_1} \times G_{2\mathbb{A}}^0/G_{2k} \\ \lambda \geq 1}} f_2(\Phi, \lambda, g_{\mathrm{st},2}^0, s, s_1) d^\times\lambda dg_{\mathrm{st},2}^0,$$

(3) $$T_{\mathrm{st},2}^1(\Phi, \omega, s_1) = \int_{(\mathbb{A}_{k_1}^1/k_1^\times)^2 \times \mathbb{A}_{k_1} \times G_{2\mathbb{A}}^0/G_{2k}} f_2(\Phi, 1, g_{\mathrm{st},2}^0, s, s_1) dg_{\mathrm{st},2}^0.$$

**Proposition (5.5)** *The integral* (5.4)(1) *converges absolutely and locally uniformly for* $\mathrm{Re}(s) > 6,\ \mathrm{Re}(s) + 4\mathrm{Re}(s_1) > 6,$ *and the integrals* (5.4)(2), (5.4)(3) *are entire functions.*

*Proof.* Let $b^0_{\mathrm{st},2+} = (n(u_0), a(\underline{\lambda}_2^{-1}, \underline{\lambda}_2))$, and $db^0_{\mathrm{st},2+} = d^\times \lambda_2 du_0$. Let $C \subset \widehat{T}^0_{1\mathbb{A}}$ be a compact set which surjects to $\widehat{T}^0_{1\mathbb{A}}/T_{1k}$. We choose $0 \le \Psi \in \mathscr{S}(V_\mathbb{A})$ so that

$$\Phi((a(t_{11}, t_{12}), g_2)x) \ll \Psi((1, a(\underline{\lambda}_2^{-1}, \underline{\lambda}_2))x)$$

for $a(t_{11}, t_{12}) \in C,\ g_2 \in \widehat{\Omega}_2 T^0_{2\eta+}$. Then

$$\Phi(\underline{\lambda} g^0_{\mathrm{st},2} w_{x,y}) \ll \Psi(\underline{\lambda}(n(t_{11}t_{12}^{-1}u_0), a(\underline{\lambda}_2^{-1}, \underline{\lambda}_2))w_{x,y}).$$

When we integrate with respect to $u_0 \in \mathbb{A}_{k_1}$, we can ignore $t_{11}t_{12}^{-1}$. So we only have to consider the integral

$$\int_{\mathbb{R}_+ \times \mathbb{A}_{k_1} \times T^0_{2\eta+}} \lambda^\sigma \alpha_{k_1}(u_0)^{\sigma_1} \Theta_{Z'^{\mathrm{ss}}_{0,2k}}(\Psi, \underline{\lambda} b^0_{\mathrm{st},2+}) d^\times \lambda db^0_{\mathrm{st},2+}. \tag{5.6}$$

If $z = \underline{\lambda} b^0_{\mathrm{st},2+} w_{x,y} = (z_{ij})$,

$$\begin{aligned} &z_{10} = z_{20} = 0, \\ &z_{11} = \underline{\lambda}\underline{\lambda}_2^{-1} x,\ z_{12} = \mathrm{Tr}_{k_1/k}(\underline{\lambda}\underline{\lambda}_2^{-1} x u_0), \\ &z_{21} = \underline{\lambda}\underline{\lambda}_2 y,\ z_{22} = \mathrm{Tr}_{k_1/k}(\underline{\lambda}\underline{\lambda}_2 y u_0). \end{aligned}$$

Let $\Psi_1 = R_{0,2}\Psi$. Then

$$\begin{aligned} &\Theta_{Z'^{\mathrm{ss}}_{0,2k}}(\Psi, \underline{\lambda} b^0_{\mathrm{st},2+}) \\ &= \sum_{x,y} \Psi_1(\underline{\lambda}\underline{\lambda}_2^{-1} x, \underline{\lambda}\underline{\lambda}_2 y, \mathrm{Tr}_{k_1/k}(\underline{\lambda}\underline{\lambda}_2^{-1} x u_0), \mathrm{Tr}_{k_1/k}(\underline{\lambda}\underline{\lambda}_2 y u_0)), \end{aligned} \tag{5.8}$$

where the sum is over all $x, y \in k_1$ such that $x, y$ are linearly independent over $k$. Let $u_3 = \mathrm{Tr}_{k_1/k}(\underline{\lambda}\underline{\lambda}_2^{-1} x u_0),\ u_4 = \mathrm{Tr}_{k_1/k}(\underline{\lambda}\underline{\lambda}_2 y u_0))$.

Suppose

$$x = x_1 + x_2\alpha_0,\ y = y_1 + y_2\alpha_0,\ u_0 = u_{01} + u_{02}\alpha_0$$

$(x_1, x_2, y_1, y_2 \in k,\ u_{01}, u_{02} \in \mathbb{A})$.

Then $du_0 = \lambda^{-2} du_3 du_4$ and

$$\underline{\lambda}^{-1} u_3 = 2\underline{\lambda}_2^{-1}(x_1 u_{01} + x_2 u_{02}\alpha_0^2),\ \underline{\lambda}^{-1} u_4 = 2\underline{\lambda}_2(y_1 u_{01} + y_2 u_{02}\alpha_0^2). \tag{5.9}$$

So

$$\frac{1}{2}\begin{pmatrix} \underline{\lambda}^{-1}\underline{\lambda}_2 u_3 \\ \underline{\lambda}^{-1}\underline{\lambda}_2^{-1} u_4 \end{pmatrix} = \begin{pmatrix} x_1 & x_2 \\ y_1 & y_2 \end{pmatrix} \begin{pmatrix} u_{01} \\ u_{02}\alpha_0^2 \end{pmatrix}.$$

Let $\Delta(x, y) = 2(x_1 y_2 - x_2 y_1)$. Since $\Delta(x, y) \ne 0$,

$$\begin{aligned} \begin{pmatrix} u_{01} \\ u_{02}\alpha_0^2 \end{pmatrix} &= \frac{\underline{\lambda}^{-1}}{\Delta(x, y)} \begin{pmatrix} y_2 & -x_2 \\ -y_1 & x_1 \end{pmatrix} \begin{pmatrix} \underline{\lambda}_2 u_3 \\ \underline{\lambda}_2^{-1} u_4 \end{pmatrix} \\ &= \frac{\underline{\lambda}^{-1}}{\Delta(x, y)} \begin{pmatrix} \underline{\lambda}_2 y_2 u_3 - \underline{\lambda}_2^{-1} x_2 u_4 \\ -\underline{\lambda}_2 y_1 u_3 + \underline{\lambda}_2^{-1} x_1 u_4 \end{pmatrix}. \end{aligned}$$

Therefore,

$$u_0 = \frac{\underline{\lambda}^{-2}}{\Delta(x,y)}((\underline{\lambda\lambda_2} y_2 u_3 - \underline{\lambda\lambda_2^{-1}} x_2 u_4) + (-\underline{\lambda\lambda_2} y_1 u_3 + \underline{\lambda\lambda_2^{-1}} x_1 u_4)\alpha_0^{-1}).$$

Let

$$f(x,y,u_3,u_4) = (y_2 u_3 - x_2 u_4) + (-y_1 u_3 + \ x_1 u_4)\alpha_0^{-1}$$

for $x = x_1 + x_2\alpha_0, y = y_1 + y_2\alpha_0 \in \mathbb{A}_{k_1}$, and $u_3, u_4 \in \mathbb{A}$. Then

$$u_0 = \frac{\underline{\lambda}^{-2}}{\Delta(x,y)} f(\underline{\lambda\lambda_2^{-1}} x, \underline{\lambda\lambda_2} y, u_3, u_4). \tag{5.10}$$

Let $L \subset k_1$ be a lattice and $F \subset \mathbb{A}_{k_1 f}$ a compact set such that

$$\Psi_1(\underline{\lambda\lambda_2^{-1}} x, \underline{\lambda\lambda_2} y, u_3, u_4) = 0$$

unless $x, y \in L$ and the finite parts of $u_3, u_4$ belong to $F$. Then there is another compact set $F_1 \subset \mathbb{A}_{k_1 f}$ such that the finite part of $f(\underline{\lambda\lambda_2^{-1}} x, \underline{\lambda\lambda_2} y, u_3, u_4)$ belongs to $F_1$ as long as the above condition is satisfied. So if $\sigma_1 \le 0$, by the consideration in [5, p. 117],

$$\alpha_{k_1}(u_0) \ll |\underline{\lambda}|_{k_1}^{2\sigma_1} \prod_{v \in \mathfrak{M}_{k_1 \infty}} \alpha_v(\underline{\lambda}^2 \Delta(x,y))^{\sigma_1} \alpha_v(f(\underline{\lambda\lambda_2^{-1}} x, \underline{\lambda\lambda_2} y, u_3, u_4))^{\sigma_1}. \tag{5.11}$$

Since each factor on the right hand side of (5.11) is a polynomial growth function of the infinite parts of $\underline{\lambda\lambda_2^{-1}} x$, $\underline{\lambda\lambda_2} y$, $u_3$, $u_4$, there exists $0 \le \Psi_2 \in \mathscr{S}(Z_{0,2\mathbb{A}})$ such that

$$\alpha(u_0)_{k_1}^{\sigma_1} \Psi_1(\underline{\lambda\lambda_2^{-1}} x, \underline{\lambda\lambda_2} y, u_3, u_4) \ll \lambda^{4\sigma_1} \Psi_2(\underline{\lambda\lambda_2^{-1}} x, \underline{\lambda\lambda_2} y, u_3, u_4)$$

as long as $\sigma_1$ ranges over a compact set.

If $\sigma_1 \ge 0$, $\alpha_{k_1}(u_0)^{\sigma_1} \le 1$. Note that if $x, y$ are linearly independent over $k$, $x, y \ne 0$. So if we define

$$\Psi_3(x_{11}, x_{21}) = \int_{\mathbb{A}^2} \Psi_2(x_{11}, x_{21}, u_3, u_4) du_3 du_4,$$

(5.6) is bounded by a constant multiple of

$$\int_{\mathbb{R}_+ \times T^0_{2\eta+}} \lambda^{\sigma-2} \sup(1, \lambda^{4\sigma_1}) \sum_{x,y \in k_1^\times} \Psi_3(\underline{\lambda\lambda_2^{-1}} x, \underline{\lambda\lambda_2} y) d^\times \lambda d^\times \lambda_2. \tag{5.12}$$

It is easy to see that (5.12) converges absolutely if $\sigma > 6$, $\sigma + 4\sigma_1 > 6$.

If $\lambda \ge 1$ or $\lambda = 1$, it is easy to see that the integral analogous to (5.12) always converges. This proves Proposition (5.5).

Q.E.D.

**Proposition (5.13)** *Let $\tau = (\nu, \nu)$. Then,*

$$\Xi_{\mathrm{st},2}(\Phi, \omega, w) \sim \frac{\varrho}{2\pi\sqrt{-1}} \int_{\mathrm{Re}(s_{\tau 1})=r_1>1} \frac{T^1_{\mathrm{st},2}(\Phi, \omega, -\frac{s_{\tau 1}-1}{2})}{s_{\tau 1}-1} \phi_{k_1}(s_{\tau 1}) \widetilde{\Lambda}(w; s_{\tau 1}, 1) ds_{\tau 1}.$$

*Proof.* We first prove that $\mathscr{E}(g_{\mathrm{st},2}, w)$ can be replaced by $G(\lambda_1, w)$ (in (2.29)). Note that the function $G(\lambda_1, w)$ does not depend on $g_2$ and therefore, is a well defined function on $X_2/T_2 \times G_{2k}$. By Proposition (2.30), $\mathscr{E}(g_{\mathrm{st},2}, w) - G(\lambda_1, w)$ is bounded by a constant multiple of a finite sum of functions of the form $\lambda_1^{c_1}\lambda_2^{c_2}$ where $c_1 > 0$ for $g_2 \in \widehat{\Omega}_2 T^0_{2\eta+}$.

Let $b_{\mathrm{st},2+} = (n(u_0), 1) d(\lambda_1, \lambda_2)$, and $db_{\mathrm{st},2+} = d^\times\lambda_1 d^\times\lambda_2 du_0$. By the argument of Proposition (5.5), we only have to prove the convergence of integrals of the form

$$\int_{\substack{\mathbb{R}^2_+ \times \mathbb{A}_{k_1} \\ \lambda_1 \le \alpha_{k_1}(u_0)^{-\frac{1}{2}}}} \lambda_1^{c_1}\lambda_2^{c_2} \Theta_{Z'^{\mathrm{ss}}_{0,2k}}(\Psi, b_{\mathrm{st},2+}) db_{\mathrm{st},2+}, \tag{5.14}$$

where $\Psi \in \mathscr{S}(V_{\mathbb{A}})$.

Since $c_1 > 0$, (5.14) is equal to

$$\frac{1}{c_1} \int_{\mathbb{R}_+ \times \mathbb{A}_{k_1}} \alpha_{k_1}(u_0)^{-\frac{c_1}{2}} \lambda_2^{c_2} \Theta_{Z'^{\mathrm{ss}}_{0,2k}}(\Psi, b^0_{\mathrm{st},2+}) db^0_{\mathrm{st},2+},$$

which converges by the proof of Proposition (5.5). Therefore,

$$\Xi_{\mathrm{st},2}(\Phi, \omega, w) \sim \int_{X_2/(T_{1k} \times G_{2k})} \Theta_{Z'^{\mathrm{ss}}_{0,2k}}(\Phi, g_{\mathrm{st},2}) G(\lambda_1, w) dg_{\mathrm{st},2}.$$

Let $\tau = (\nu, \nu)$. Note that $\Theta_{Z'^{\mathrm{ss}}_{0,2k}}(\Phi, g_{\mathrm{st},2}) = \Theta_{Z'^{\mathrm{ss}}_{0,2k}}(\Phi, g^0_{\mathrm{st},2})$. So

$$\begin{aligned}
&\int_{X_2/(T_{1k} \times G_{2k})} \omega(g^0_{\mathrm{st},2}) \lambda_1^{s_{\tau 1}-1} \alpha_{k_1}(u_0)^{s_1} \Theta_{Z'^{\mathrm{ss}}_{0,2k}}(\Phi, g^0_{\mathrm{st},2}) d^\times\lambda_1 dg^0_{\mathrm{st},2} \\
&= \int_{(\mathbb{A}^1_{k_1}/k_1^\times)^2 \times \mathbb{A}_{k_1} \times G^0_{2\mathbb{A}}/G_{2k}} \omega(g^0_{\mathrm{st},2}) \frac{\alpha_{k_1}(u_0)^{-\frac{s_{\tau 1}-1}{2}}}{s_{\tau 1}-1} \Theta_{Z'^{\mathrm{ss}}_{0,2k}}(\Phi, g^0_{\mathrm{st},2}) dg^0_{\mathrm{st},2} \\
&= \frac{T^1_{\mathrm{st},2}(\Phi, \omega, -\frac{s_{\tau 1}-1}{2})}{s_{\tau 1}-1}.
\end{aligned}$$

Therefore,

$$\begin{aligned}
&\int_{X_2/T_2 \times G_{2k}} \omega(g^0_{\mathrm{st},2}) \Theta_{Z'^{\mathrm{ss}}_{0,2k}}(\Phi, g_{\mathrm{st},2}) G(\lambda_1, w) dg_{\mathrm{st},2} \\
&\sim \frac{\varrho}{2\pi\sqrt{-1}} \int_{\mathrm{Re}(s_{\tau 1})=r_1>1} \frac{T^1_{\mathrm{st},2}(\Phi, \omega, -\frac{s_{\tau 1}-1}{2})}{s_{\tau 1}-1} \phi_{k_1}(s_{\tau 1}) \widetilde{\Lambda}(w; s_{\tau 1}, 1) ds_{\tau 1}.
\end{aligned}$$

This proves the proposition.

Q.E.D.

We define the adjusting term corresponding to $V^{\rm ss}_{{\rm st},2k}$ as follows.

**Definition (5.15)**

(1) $$T_{{\rm st},2}(\Phi,\omega,s)=\left.\frac{d}{ds_1}\right|_{s_1=0}T_{{\rm st},2}(\Phi,\omega,s,s_1),$$

(2) $$T_{{\rm st},2+}(\Phi,\omega,s)=\left.\frac{d}{ds_1}\right|_{s_1=0}T_{{\rm st},2+}(\Phi,\omega,s,s_1),$$

(3) $$T^1_{{\rm st},2}(\Phi,\omega)=\left.\frac{d}{ds_1}\right|_{s_1=0}T^1_{{\rm st},2}(\Phi,\omega,s_1).$$

It follows from Proposition (5.5) that $T_{{\rm st},2}(\Phi,\omega,s)$ is holomorphic for ${\rm Re}(s)>6$ and $T_{{\rm st},2+}(\Phi,\omega,s)$ is an entire function.

We show that $T_{{\rm st},2}(\Phi,\omega,s,s_1)$ has an Euler product. Let $\delta_{{\rm st},2}(\omega)=\delta(\omega_1(\omega_2\circ {\rm N}_{k_1/k})^{-1})$.

**Proposition (5.16)**

$$\begin{aligned}&T_{{\rm st},2}(\Phi,\omega,s,s_1)\\&=\delta_{{\rm st},2}(\omega)\int_{\mathbb{R}_+\times\mathbb{A}_{k_1}\times G^0_{2\mathbb{A}}}\omega_2(\det g_2)\lambda^s\alpha_{k_1}(u_0)^{s_1}\Phi(\underline{\lambda}g_2w_2)d^\times\lambda du_0dg_2.\end{aligned}$$

*Proof.* The group $G_{w_2k}\cong {\rm GL}(1)_{k_1}\times{\rm GL}(1)_{k_1}$ is a subgroup of $T_1\times G_{2k}$, and the imbedding is given by $(t_{11},t_{12})\to(a(t_{11},t_{12}),A(t_{12},t_{11}))$ (see §1). We use $d^\times t_{11}d^\times t_{12}$ as the measure on $G_{w_2\mathbb{A}}$. Let $G^0_{w_2\mathbb{A}}=(\mathbb{A}^1_{k_1})^2$. The volume of $G^0_{w_2\mathbb{A}}/G_{w_2k}$ is 1. It is easy to see that $(T_{1k}\times G_{2k})/G_{w_2k}\cong Z'^{\rm ss}_{0,2k}$ and this is given by $(t_1,g_2)\to(t_1,g_2)w_2$. Also it is easy to see that

$$\omega(a(t_{11},t_{12}),A(t_{12},t_{11}))=\omega_1(t_{11}t_{12})\omega_2({\rm N}_{k_1/k}(t_{11}t_{12})).$$

Since $G^0_{2\mathbb{A}}$ is isomorphic to its image in $\mathbb{R}_+\times\widehat{T}^0_{1\mathbb{A}}\times G^0_{2\mathbb{A}}/G_{w_2\mathbb{A}}$,

$$\begin{aligned}&T_{{\rm st},2}(\Phi,\omega,s,s_1)\\&=\int_{\mathbb{R}_+\times\mathbb{A}_{k_1}\times\widehat{T}^0_{1\mathbb{A}}\times G^0_{2\mathbb{A}}/G_{w_2k}}\omega(g^0_{{\rm st},2})\lambda^s\alpha_{k_1}(u_0)^{s_1}\Phi(\underline{\lambda}g^0_{{\rm st},2}w_2)d^\times\lambda dg^0_{{\rm st},2}\\&=\delta_{{\rm st},2}(\omega)\int_{\mathbb{R}_+\times\mathbb{A}_{k_1}\times G^0_{2\mathbb{A}}}\omega_2(\det g_2)\lambda^s\alpha_{k_1}(u_0)^{s_1}\Phi(\underline{\lambda}g_2w_2).\end{aligned}$$

Q.E.D.

Clearly, the above integral can be rewritten as an integral on $\mathbb{A}_{k_1}\times G_{2\mathbb{A}}$. Therefore, it has an Euler product.

## §6 Contributions from $\mathfrak{p}_1,\mathfrak{p}_3$

We consider contributions from $\mathfrak{p}_1, \mathfrak{p}_3$ in this section. We first consider $\mathfrak{p}_1$.

**Lemma (6.1)** $\widetilde{\Xi}_{\mathfrak{p}_1}(\Phi, \omega, w) \sim 0$.

*Proof.* For $\alpha = (\alpha_1, \alpha_2) \in k_1 \times k$, we define

$$f_\alpha(t^0) = \int_{N_{\mathbb{A}}/N_k} \Theta_{Y_{1k}^{\mathrm{ss}}}(\Phi, t^0 n(u)) < \alpha n(u) > du,$$

.

By the Parseval formula,

$$\widetilde{\Xi}_{\mathfrak{p}_1}(\Phi, \omega, w) = \int_{T_{\mathbb{A}}^0/T_k} \sum_{\alpha \in k_1 \times k \setminus \{(0,0)\}} f_\alpha(t^0) \mathscr{E}_\alpha(t^0, w) \lambda_1^2 \lambda_2^2 d^\times t^0.$$

**Lemma (6.2)** *For any* $M \gg 0$,

$$\sum_{\alpha \in k_1 \times k \setminus \{(0,0)\}} |f_{1,\alpha}(t^0)| \ll \mathrm{rd}_{2,M}(\lambda_1, \lambda_2).$$

*Proof.* Let

$$A(x, t^0, u) = (\gamma_{12}(t^0) x_{12}, \gamma_{21}(t^0) x_{21}, \gamma_{22}(t^0)(x_{22} + \mathrm{Tr}_{k_1/k}(x_{21} u_1) + x_{12} u_2)).$$

Since $N_{\mathbb{A}}/N_k$ is compact,

$$f_{1,\alpha}(t^0) = \sum_{\substack{x_{12} \in k^\times, \\ x_{21} \in k_1^\times}} \int_{N_{\mathbb{A}}/N_k} \sum_{x_{22} \in k} \widetilde{R}_1 \Phi(A(x, t^0, u)) < \alpha n(u) > du. \tag{6.3}$$

Let

$$u_1' = u_1'(x, u) = x_{21} u_1 \in \mathbb{A}_{k_1},\ u_2' = u_2'(x, u) = \mathrm{Tr}_{k_1/k}(x_{21} u_1) + x_{12} u_2 \in \mathbb{A}.$$

Then

$$u_1 = x_{21}^{-1} u_1',\ u_2 = x_{12}^{-1}(u_2' - \mathrm{Tr}_{k_1/k}(u_1')).$$

So

$$\mathrm{Tr}_{k_1/k}(\alpha_1 u_1) + \alpha_2 u_2 = \mathrm{Tr}_{k_1/k}((\alpha_1 x_{21}^{-1} - \alpha_2 x_{12}^{-1}) u_1') + \alpha_2 x_{12}^{-1} u_2'.$$

Therefore, for fixed $x_{12}, x_{21}$,

$$\int_{N_{\mathbb{A}}/N_k} \sum_{x_{22} \in k} \widetilde{R}_1 \Phi(A(x, t^0, u)) < \alpha n(u) > du = 0$$

unless $\alpha_1 = x_{12}^{-1} x_{21} \alpha_2$.

Suppose $\alpha_1 = x_{12}^{-1}x_{21}\alpha_2$. This implies that the condition $\alpha \neq (0,0)$ is equivalent to the condition $\alpha_1, \alpha_2 \neq 0$. Let $\Phi_1$ be the partial Fourier transform of $\widetilde{R}_1\Phi$ with respect to the coordinate $x_{22}$ and the character $< >$. Then

$$\int_{N_{\mathbb{A}}/N_k} \sum_{x_{22}\in k} \widetilde{R}_1\Phi(A(x,t^0,u)) < \alpha n(u) > du$$
$$= |\gamma_{22}(t^0)|^{-1}\Phi_1(\gamma_{12}(t^0)x_{12}, \gamma_{21}(t^0)x_{21}, \gamma_{22}(t^0)^{-1}\alpha_2 x_{12}^{-1}).$$

So

$$\sum_{\alpha\in k_1\times k\setminus\{(0,0)\}} |f_{1,\alpha}(t^0)|$$
$$\ll (\lambda_1\lambda_2)^{-1} \sum_{\substack{x_{12},\alpha_2\in k^\times,\\ x_{21}\in k_1^\times}} |\Phi_1(\gamma_{12}(t^0)x_{12}, \gamma_{21}(t^0)x_{21}, \gamma_{22}(t^0)^{-1}\alpha_2 x_{12}^{-1})|.$$

Therefore, by Lemma (1.2.6) [5], for any $N_1, N_2, N_3 \geq 1$,

$$\sum_{\alpha\in k_1\times k\setminus\{(0,0)\}} |f_{1,\alpha}(t^0)| \ll (\lambda_1\lambda_2)^{-1}(\lambda_1\lambda_2^{-1})^{-N_1}(\lambda_2)^{-2N_2}(\lambda_1\lambda_2)^{N_3}.$$

Since the convex hull of $\{(1,-1),(0,2),(-1,-1)\} \subset \mathbb{R}^2$ contains a neighborhood of the origin of $\mathbb{R}^2$, the lemma follows.

Q.E.D.

We continue the proof of (6.1). By Proposition (2.30), there exists a slowly increasing function $h(\lambda_1,\lambda_2)$ and $\delta > 0$ such that on any vertical strip in $\mathrm{Re}(w) > L(\rho) - \delta$,

$$\mathscr{E}_\alpha(t^0, w) \ll h(\lambda_1,\lambda_2)$$

and this bound is uniform with respect to $\alpha \in k_1 \times k \setminus \{(0,0)\}$. So for any $M \gg 0$,

$$\widetilde{\Xi}_{\mathfrak{p}_1}(\Phi,\omega,w) \ll \int_{\mathbb{R}_+^2} \lambda_1\lambda_2 h(\lambda_1,\lambda_2)\mathrm{rd}_{2,M}(\lambda_1,\lambda_2)d^\times\lambda_1 d^\times\lambda_2 < \infty.$$

This proves Lemma (6.1).

Q.E.D.

Consider $\Xi_{\mathfrak{p}_1}(\Phi,\omega,w)$.

For $\omega = (\omega_1,\omega_2)$ as before, we define $\delta_{\mathfrak{d}_1}(\omega) = \delta(\omega_1|_{\mathbb{A}^\times}\omega_2^{-1})\delta(\omega_1^2)$ and $\omega_{\mathfrak{d}_1} = (\omega_2,\omega_1)$. For each Weyl group element $\tau$, we define

$$\text{(6.4)}\qquad \Sigma_{\mathfrak{d}_1}(\Phi,\omega,s_\tau) = \frac{\delta_{\mathfrak{d}_1}(\omega)}{2}\Sigma_{1,1}\left(R_1\Phi, \omega_{\mathfrak{d}_1}, s_{\tau 1}, \frac{s_{\tau 1}+s_{\tau 2}}{2}\right).$$

**Proposition (6.5)** *Let $\tau = (\nu,\nu)$. Then,*

$$\Xi_{\mathfrak{p}_1}(\Phi,\omega,w) \sim \frac{\varrho}{2\pi\sqrt{-1}} \int_{\mathrm{Re}(s_{\tau 1})=r_1>1} \Sigma_{\mathfrak{d}_1}(\Phi,\omega,s_{\tau 1},1)\phi_{k_1}(s_{\tau 1})\widetilde{\Lambda}(w;s_{\tau 1},1)ds_{\tau 1}.$$

*Proof.* Note that

$$\gamma_{12}(t^0) = \underline{\lambda}_1\underline{\lambda}_2^{-1}\mathrm{N}_{k_1/k}(t_{12})t_{21},\ \gamma_{21}(t^0) = \underline{\lambda}_2 t_{11}t_{12}^{\sigma}t_{22}.$$

If $f(q_1, q_2)$ is a function of $(q_1, q_2) \in \mathbb{A}^1/k^\times \times \mathbb{A}^1_{k_1}/k_1^\times$,

$$\begin{aligned}&\int_{\widehat{T}^0_{\mathbb{A}}/T_k} \omega(\widehat{t}^0)f(\mathrm{N}_{k_1/k}(t_{12})t_{21}, t_{11}t_{12}^{\sigma}t_{22})d^\times\widehat{t}^0\\ &= \delta_{\mathfrak{d}_1}(\omega)\int_{\mathbb{A}^1/k^\times\times\mathbb{A}^1_{k_1}/k_1^\times} \omega_2(q_1)\omega_1(q_2)f(q_1,q_2)d^\times q_1 d^\times q_2\end{aligned}.$$

Let $t' = (t'_1, t'_2) = (\underline{\mu}_1 q_1, \underline{\mu}_{2k_1} q_2) \in \mathbb{A}^\times \times \mathbb{A}^\times_{k_1}$. We make the change of variables $\mu_1 = \lambda_1\lambda_2^{-1},\ \mu_2 = \lambda_2^2$. Then $d^\times\lambda_1 d^\times\lambda_2 = \frac{1}{2}d^\times\mu_1 d^\times\mu_2$. Also for $s_\tau$ as in (2.26),

$$\lambda_1^{s_{\tau 1}-1}\lambda_2^{s_{\tau 2}-1}(t^0)^{-2\rho}|\gamma_{22}(t^0)|^{-1} = \mu_1^{s_{\tau 1}}\mu_2^{\frac{s_{\tau 1}+s_{\tau 2}}{2}}.$$

Therefore, if $\mathrm{Re}(s_{\tau 1}), \mathrm{Re}(s_{\tau 2}) > 1$,

$$\begin{aligned}&\int_{T^0_{\mathbb{A}}/T_k} \omega(t^0)\Theta_{1,1}(R_1\Phi, \gamma_{12}(t^0), \gamma_{21}(t^0))|\gamma_{22}(t^0)|^{-1}\lambda_1^{s_{\tau 1}-1}\lambda_2^{s_{\tau 2}-1}d^\times t^0\\ &= \frac{\delta_{\mathfrak{d}_1}(\omega)}{2}\int_{\mathbb{A}^\times/k^\times\times\mathbb{A}^\times_{k_1}/k_1^\times} \omega_2(t'_1)\omega_1(t'_2)\Theta_{1,1}(R_1\Phi, t'_1, t'_2)|t'_1|^{s_{\tau 1}}|t'_2|^{\frac{s_{\tau 1}+s_{\tau 2}}{2}}d^\times t'_1 d^\times t'_2\\ &= \Sigma_{\mathfrak{d}_1}(\Phi, \omega, s_\tau).\end{aligned}$$

So

$$\Xi_{\mathfrak{p}_1}(\Phi, \omega, w) = \sum_\tau \left(\frac{1}{2\pi\sqrt{-1}}\right)^2 \int_{\substack{\mathrm{Re}(s_\tau)=(r_1,r_2)\\ r_1,r_2>1}} \Sigma_{\mathfrak{d}_1}(\Phi, \omega, s_\tau)\widetilde{\Lambda}_\tau(w; s_\tau)ds_\tau.$$

If $\tau \neq (\nu, \nu)$, we can choose $r_1, r_2$ close enough to 1 so that $\widetilde{L}(r) < L(\rho)$. So

$$\left(\frac{1}{2\pi\sqrt{-1}}\right)^2 \int_{\substack{\mathrm{Re}(s_\tau)=(r_1,r_2)\\ r_1,r_2>1}} \Sigma_{\mathfrak{d}_1}(\Phi, \omega, s_\tau)\widetilde{\Lambda}_\tau(w; s_\tau)ds_\tau \sim 0.$$

Let $\tau = (\nu, \nu)$. Since $C > 4$, $\widetilde{L}(2, \frac{1}{2}) = 2 + \frac{1}{2}C < 1 + C = L(\rho)$. So

$$\begin{aligned}&\left(\frac{1}{2\pi\sqrt{-1}}\right)^2 \int_{\substack{\mathrm{Re}(s_\tau)=(r_1,r_2)\\ r_1,r_2>1}} \Sigma_{\mathfrak{d}_1}(\Phi, \omega, s_\tau)\widetilde{\Lambda}_\tau(w; s_\tau)ds_\tau\\ &= \left(\frac{1}{2\pi\sqrt{-1}}\right)^2 \int_{\mathrm{Re}(s_\tau)=(2,\frac{1}{2})} \Sigma_{\mathfrak{d}_1}(\Phi, \omega, s_\tau)\widetilde{\Lambda}_\tau(w; s_\tau)ds_\tau\\ &\quad + \frac{\varrho}{2\pi\sqrt{-1}}\int_{\mathrm{Re}(s_{\tau 1})=r_1>1} \Sigma_{\mathfrak{d}_1}(\Phi, \omega, s_{\tau 1}, 1)\phi_{k_1}(s_{\tau 1})\widetilde{\Lambda}(w; s_{\tau 1}, 1)ds_{\tau 1}\\ &\sim \frac{\varrho}{2\pi\sqrt{-1}}\int_{\mathrm{Re}(s_{\tau 1})=r_1>1} \Sigma_{\mathfrak{d}_1}(\Phi, \omega, s_{\tau 1}, 1)\phi_{k_1}(s_{\tau 1})\widetilde{\Lambda}(w; s_{\tau 1}, 1)ds_{\tau 1}.\end{aligned}$$

This proves the proposition.

Q.E.D.

Next we consider $\mathfrak{p}_3$.

**Proposition (6.6)**

$$\Xi_{\mathfrak{p}_3}(\Phi,\omega,w)=\sum_\tau \frac{\delta_\#(\omega)}{2\pi\sqrt{-1}}\int_{\mathrm{Re}(s_{\tau 1})=r>1}\Sigma_1(R_3\Phi,s_{\tau 1}+1)\widetilde{\Lambda}_\tau(w;s_{\tau 1},s_{\tau 1})ds_{\tau 1}.$$

*Proof.* If $f(q)$ is a function of $q\in\mathbb{A}^1/k^\times$,

$$\int_{\widehat{T}^0_{\mathbb{A}}/T_k} f(\mathrm{N}_{k_1/k}(t_{12})t_{22})d^\times\widehat{t}^0=\delta_\#(\omega)\int_{\mathbb{A}^1/k^\times} f(q)d^\times q. \tag{6.7}$$

For $s_\tau$ as in (2.26),

$$\lambda_1^{s_{\tau 1}-1}\lambda_2^{s_{\tau 2}-1}(t^0)^{-2\rho}=\lambda_1^{s_{\tau 1}+1}\lambda_2^{s_{\tau 2}+1}.$$

Since $|\gamma_{22}(t^0)|=\lambda_1\lambda_2$, we make the change of variables $\mu_1=\lambda_1\lambda_2,\ \mu_2=\lambda_2$. Then $\lambda_1^{s_{\tau 1}+1}\lambda_2^{s_{\tau 2}+1}=\mu_1^{s_{\tau 1}+1}\mu_2^{s_{\tau 2}-s_{\tau 1}}$. Since $\gamma_{22}(t^0)$ does not depend on $\mu_2$, by the Mellin inversion formula (or by using (3.5.17) [5] more precisely),

$$\Xi_{\mathfrak{p}_3}(\Phi,\omega,w)=\sum_\tau \frac{\delta_\#(\omega)}{2\pi\sqrt{-1}}\int_{\mathrm{Re}(s_{\tau 1})=r>1}\Sigma_1(R_3\Phi,s_{\tau 1}+1)\widetilde{\Lambda}_\tau(w;s_{\tau 1},s_{\tau 1})ds_{\tau 1}.$$

Q.E.D.

It turns out that $\Xi_{\mathfrak{p}_3}(\Phi,\omega,w)$ will cancelled out with a distribution which arises when we move the contour defining $\Xi_{\mathfrak{p}_{41}}(\Phi,\omega,w)$.

## §7 Contributions from $\mathfrak{p}_2,\mathfrak{p}_{41},\mathfrak{p}_{42}$

We consider contributions from $\mathfrak{p}_2,\mathfrak{p}_{41},\mathfrak{p}_{42}$ in this section. We first have to define certain distributions.

We identify $Z_2$ with $W$. We denote the Fourier transform which we defined in §3 by $\mathscr{F}_W$.

We first consider $\Xi_{\mathfrak{p}_2+}(\Phi,\omega,w)$, $\widehat{\Xi}_{\mathfrak{p}_2+}(\Phi,\omega,w)$.

If $f$ is a function of $q\in\mathbb{A}^1/k^\times$,

$$\int_{(\mathbb{A}^1/k^\times)^2}\omega_2(t_{21}t_{22})f(t_{22})d^\times t_{21}d^\times t_{22}=\delta(\omega_2)\int_{\mathbb{A}^1/k^\times} f(q)d^\times q. \tag{7.1}$$

Let $g_1\in G^0_{1\mathbb{A}},\ t_2=a(\underline{\lambda}_2^{-1}t_{21},\underline{\lambda}_2 t_{22})\in T^0_{2\mathbb{A}}$. We define $d^\times t_2=d^\times\lambda_2 d^\times t_{21}d^\times t_{22}$. Then it is easy to see that

$$\begin{aligned}&\sigma_{\mathfrak{p}_2}(g_1,t_2)=\lambda_2^2,\ \kappa_{\mathfrak{d}1}(g_1,t_2)=\lambda_2^{-4},\ \sigma_{\mathfrak{p}_2}(g_1,t_2)\kappa_{\mathfrak{d}1}(g_1,t_2)=\lambda_2^{-2},\\ &\kappa_{\mathfrak{d}2}(g_1,t_2)=1,\ \chi(g_1,t_2)=1\end{aligned} \tag{7.2}$$

(the last equation is because $W$ is an irreducible representation). For the definition of these values, see [5, pp. 73, 86].

**Proposition (7.3)** (1) $\Xi_{\mathfrak{p}_2+}(\Phi,\omega,w) \sim C_G\Lambda(w;\rho)\delta(\omega_2)Z_{W+}(R_2\Phi,(1,\omega_1),2)$.
(2) $\widehat{\Xi}_{\mathfrak{p}_2+}(\Phi,\omega,w) \sim C_G\Lambda(w;\rho)\delta(\omega_2)Z_{W+}(\mathscr{F}_W R_2\Phi,(1,\omega_1^{-1}),2)$.

*Proof.* Let

$$A(\Phi,g_1,\lambda_2,q) = \omega_1(\det g_1)\lambda_2^2\Theta_{Z_{2k}^{\mathrm{ss}}}(R_2\Phi,\underline{\lambda}_2 q g_1)\mathscr{E}_{N_2}((g_1,a(\underline{\lambda}_2^{-1},\underline{\lambda}_2),w),$$

$$B(\Phi,g_1,\lambda_2,q) = \omega_1(\det g_1)\lambda_2^{-2}\Theta_{Z_{2k}^{\mathrm{ss}}}(\mathscr{F}_W R_2\Phi,\underline{\lambda}_2^{-1}q^{-1}g_1^{\iota})\mathscr{E}_{N_2}((g_1,a(\underline{\lambda}_2^{-1},\underline{\lambda}_2),w).$$

Then by (2.19), (7.1), (7.2),

$$\Xi_{\mathfrak{p}_2+}(\Phi,\omega,w) = \delta(\omega_2)\int_{\substack{G_{1\mathbb{A}}^0/G_{1k}\times\mathbb{R}_+\times\mathbb{A}^1/k^\times \\ \lambda_2\geq 1}} A(\Phi,g_1,\lambda_2,q)dg_1 d^\times\lambda_2 d^\times q,$$

$$\widehat{\Xi}_{\mathfrak{p}_2+}(\Phi,\omega,w) = \delta(\omega_2)\int_{\substack{G_{1\mathbb{A}}^0/G_{1k}\times\mathbb{R}_+\times\mathbb{A}^1/k^\times \\ \lambda_2\leq 1}} B(\Phi,g_1,\lambda_2,q)dg_1 d^\times\lambda_2 d^\times q.$$

Roughly speaking, (3.4.31)(1) [5] says that if the theta series has a bound by power functions such that the integral of the above kind without the Eisenstein series, replacing the theta series by the bound, converges absolutely, then we can replace the Eisenstein series by $C_G\Lambda(w;\rho)$. This is the case because of (3.5). Therefore, the proposition follows.

Q.E.D.

Next, we consider $\Xi_{\mathfrak{p}_2\#}(\Phi,\omega,w)$, $\widehat{\Xi}_{\mathfrak{p}_2\#}(\Phi,\omega,w)$.

**Proposition (7.4)** (1) $\Xi_{\mathfrak{p}_2\#}(\Phi,\omega,w) \sim \frac{\varrho\delta_\#(\omega)}{2}\Lambda(w;\rho)R_2\Phi(0)$
(2) $\widehat{\Xi}_{\mathfrak{p}_2\#}(\Phi,\omega,w) \sim -\frac{\varrho\delta_\#(\omega)}{2}\Lambda(w;\rho)\mathscr{F}_W R_2\Phi(0)$

*Proof.* By (7.2), and the Mellin inversion formula,

$$\begin{aligned}\Xi_{\mathfrak{p}_2\#}(\Phi,\omega,w) &= \delta_\#(\omega)R_2\Phi(0)\int_{\substack{G_{1\mathbb{A}}^0/G_{1k}\times T_{2\mathbb{A}}^0/T_{2k} \\ \lambda_2\leq 1}} \lambda_2^2\mathscr{E}_{N_2}((g_1,t_2),w)dg_1 d^\times t_2 \\ &= \sum_{\tau=(1,\tau_2)}\frac{\delta_\#(\omega)R_2\Phi(0)}{2\pi\sqrt{-1}}\int_{\mathrm{Re}(s_{\tau 2})=r_2>1}\frac{\widetilde{\Lambda}_{(1,\tau_2)}(w;-1,s_{\tau 2})}{s_{\tau 2}+1}ds_{\tau 2},\end{aligned}$$

$$\begin{aligned}\widehat{\Xi}_{\mathfrak{p}_2\#}(\Phi,\omega,w) &= \delta_\#(\omega)\mathscr{F}_W R_2\Phi(0)\int_{\substack{G_{1\mathbb{A}}^0/G_{1k}\times T_{2\mathbb{A}}^0/T_{2k} \\ \lambda_2\leq 1}} \lambda_2^{-2}\mathscr{E}_{N_2}((g_1,t_2),w)dg_1 d^\times t_2 \\ &= \sum_{\tau=(1,\tau_2)}\frac{\delta_\#(\omega)R_2\Phi(0)}{2\pi\sqrt{-1}}\int_{\mathrm{Re}(s_{\tau 2})=r_2>3}\frac{\widetilde{\Lambda}_{(1,\tau_2)}(w;-1,s_{\tau 2})}{s_{\tau 2}-3}ds_{\tau 2}.\end{aligned}$$

If $\tau_2 = 1$, we can choose $r_2 \gg 0$ and ignore the corresponding term. So we consider $\tau_2 = \nu$. Then

$$\begin{aligned}\Xi_{\mathfrak{p}_2\#}(\Phi,\omega,w) &\sim \frac{\delta_\#(\omega)R_2\Phi(0)}{2\pi\sqrt{-1}}\int_{\mathrm{Re}(s_{\tau 2})=\frac{1}{2}}\frac{\widetilde{\Lambda}_{(1,\tau_2)}(w;-1,s_{\tau 2})}{s_{\tau 2}+1}ds_{\tau 2} \\ &\quad + \frac{\varrho\delta_\#(\omega)}{2}\Lambda(w;\rho)R_2\Phi(0) \\ &\sim \frac{\varrho\delta_\#(\omega)}{2}\Lambda(w;\rho)R_2\Phi(0)\end{aligned}$$

Since $\rho$ does not satisfy the condition $s_{\tau 2}=3$, we can use the passing principle (3.6.1) [5], and we don't have to worry about the denominator $s_{\tau 2}-3$ for $\widehat{\Xi}_{\mathfrak{p}_2\#}(\Phi,\omega,w)$. So

$$\begin{aligned}\widehat{\Xi}_{\mathfrak{p}_2\#}(\Phi,\omega,w) &\sim \frac{\delta_{\#}(\omega)\mathscr{F}_W R_2\Phi(0)}{2\pi\sqrt{-1}}\int_{\mathrm{Re}(s_{\tau 2})=\frac{1}{2}} \frac{\widetilde{\Lambda}_{(1,\tau_2)}(w;-1,s_{\tau 2})}{s_{\tau 2}-3}ds_{\tau 2} \\ &\quad -\frac{\varrho\delta_{\#}(\omega)}{2}\Lambda(w;\rho)\mathscr{F}_W R_2\Phi(0) \\ &\sim -\frac{\varrho\delta_{\#}(\omega)}{2}\Lambda(w;\rho)\mathscr{F}_W R_2\Phi(0).\end{aligned}$$

This proves the proposition.

Q.E.D.

Next, we consider $\Xi_{\mathfrak{p}_{41}}(\Phi,\omega,w),\Xi_{\mathfrak{p}_{42}}(\Phi,\omega,w)$ (see [5, p. 71] for the definition of $\Phi_{\mathfrak{p}}$).

**Proposition (7.5)**

(1) $\Xi_{\mathfrak{p}_{41}}(\Phi,\omega,w)$

$$\sim \delta_{\#}(\omega)\sum_{\tau}\left(\frac{1}{2\pi\sqrt{-1}}\right)^2\int_{\substack{\mathrm{Re}(s_\tau)=(r_1,r_2)\\ r_2>r_1>1}} \frac{\Sigma_1(R_4\Phi_{\mathfrak{p}_{41}},s_{\tau 1}+1)}{s_{\tau 2}-s_{\tau 1}}\widetilde{\Lambda}_\tau(w;s_\tau)ds_\tau,$$

(2) $\Xi_{\mathfrak{p}_{42}}(\Phi,\omega,w)$

$$\sim \frac{\varrho\delta_{\#}(\omega)}{2\pi\sqrt{-1}}\int_{\mathrm{Re}(s_{\tau 1})=r_1>1}\frac{\Sigma_1(R_4\Phi_{\mathfrak{p}_{42}},s_{\tau 1}+1)}{s_{\tau 1}-1}\phi_{k_1}(s_{\tau 1})\widetilde{\Lambda}(w;s_{\tau 1},1)ds_{\tau 1},$$

*where* $\tau=(\nu,\nu)$ *in* (2).

*Proof.* Let $\tau$ be a Weyl group element. By easy computations,

$$\begin{aligned}&\sigma_{\mathfrak{p}_{41}}(t^0)=\lambda_1^2\lambda_2^2, \\ &\sigma_{\mathfrak{p}_{42}}(t^0)=(\lambda_2^2\lambda_2^{-4})^{-1}\lambda_1^2=\lambda_1^2\lambda_2^2, \\ &(t^0)^{\tau z+\rho}=\lambda_1^{s_{\tau 1}-1}\lambda_2^{s_{\tau 2}-1}, \\ &\theta_{\mathfrak{d}_2}(t^0)^{\tau z+\rho}=\lambda_1^{s_{\tau 1}-1}(\lambda_2^{-1})^{s_{\tau 2}-1}=\lambda_1^{s_{\tau 1}-1}\lambda_2^{1-s_{\tau 2}}.\end{aligned}\tag{7.6}$$

Let $\mu=\lambda_1\lambda_2,\ \lambda_2=\lambda_2$. Then

$$\begin{aligned}&\sigma_{\mathfrak{p}_{41}}(t^0)(t^0)^{\tau z+\rho}=\lambda_1^{s_{\tau 1}+1}\lambda_2^{s_{\tau 2}+1}=\mu^{s_{\tau 1}+1}\lambda_2^{s_{\tau 2}-s_{\tau 1}}, \\ &\sigma_{\mathfrak{p}_{42}}(t^0)\theta_{\mathfrak{d}_2}(t^0)^{\tau z+\rho}=\lambda_1^{s_{\tau 1}+1}\lambda_2^{3-s_{\tau 2}}=\mu^{s_{\tau 1}+1}\lambda_2^{2-s_{\tau 1}-s_{\tau 2}}.\end{aligned}\tag{7.7}$$

Since

$$\int_{\substack{\mathbb{R}_+\times\mathbb{A}^1/k^\times\\ \lambda_2\le 1}}\mu^{s_{\tau 1}+1}\lambda_2^{s_{\tau 2}-s_{\tau 1}}\Theta_{Z_4}(R_4\Phi_{\mathfrak{p}_{41}},\underline{\mu}q)d^\times\mu d^\times q=\frac{\Sigma_1(R_4\Phi_{\mathfrak{p}_{41}},s_{\tau 1}+1)}{s_{\tau 2}-s_{\tau 1}},$$

$$\int_{\substack{\mathbb{R}_+\times\mathbb{A}^1/k^\times\\ \lambda_2\ge 1}}\mu^{s_{\tau 1}+1}\lambda_2^{2-s_{\tau 1}-s_{\tau 2}}\Theta_{Z_4}(R_4\Phi_{\mathfrak{p}_{42}},\underline{\mu}q)d^\times\mu d^\times q=\frac{\Sigma_1(R_4\Phi_{\mathfrak{p}_{42}},s_{\tau 1}+1)}{s_{\tau 1}+s_{\tau 2}-2},$$

by (2.21), (6.7), (7.6), (7.7),

$$\begin{aligned}
&\Xi_{\mathfrak{p}_{41}}(\Phi,\omega,w)\\
&= \delta_{\#}(\omega)\sum_{\tau}\left(\frac{1}{2\pi\sqrt{-1}}\right)^2 \int_{\substack{\mathrm{Re}(s_\tau)=(r_1,r_2)\\ r_2>r_1>1}} \frac{\Sigma_1(R_4\Phi_{\mathfrak{p}_{41}},s_{\tau 1}+1)}{s_{\tau 2}-s_{\tau 1}}\widetilde{\Lambda}_\tau(w;s_\tau)ds_\tau,\\
&\Xi_{\mathfrak{p}_{42}}(\Phi,\omega,w)\\
&= \delta_{\#}(\omega)\sum_{\tau}\left(\frac{1}{2\pi\sqrt{-1}}\right)^2 \int_{\substack{\mathrm{Re}(s_\tau)=(r_1,r_2)\\ r_1,r_2>1}} \frac{\Sigma_1(R_4\Phi_{\mathfrak{p}_{42}},s_{\tau 1}+1)}{s_{\tau 1}+s_{\tau 2}-2}\widetilde{\Lambda}_\tau(w;s_\tau)ds_\tau.
\end{aligned}$$

Consider $\Xi_{\mathfrak{p}_{42}}(\Phi,\omega,w)$. If $\tau\neq(\nu,\nu)$, by choosing $r_1\gg 0$ or $r_2\gg 0$, we can ignore the corresponding term. Let $\tau=(\nu,\nu)$. Then by moving the contour from $r_1,r_2>1$ to $r_1=2,\ r_2=\frac{1}{2}$, we get (2). This proves the proposition.

Q.E.D.

Note that by (3.6), (7.3), (7.4),

$$\text{(7.6)}\qquad \begin{aligned}
&\Xi_{\mathfrak{p}_2+}(\Phi,\omega,w)+\widehat{\Xi}_{\mathfrak{p}_2+}(\Phi,\omega,w)+\widehat{\Xi}_{\mathfrak{p}_2\#}(\Phi,\omega,w)-\Xi_{\mathfrak{p}_2\#}(\Phi,\omega,w)\\
&\sim C_G\Lambda(w;\rho)\delta(\omega_2)Z_{W,(0)}(R_2\Phi,(1,\omega_1),2).
\end{aligned}$$

## §8 The principal part formula

**Proposition (8.1)** *Let* $\tau=(\nu,\nu)$. *Then,*

$$\begin{aligned}
&\Xi_{\mathfrak{p}_3}(\Phi,\omega,w)-\Xi_{\mathfrak{p}_{41}}(\Phi,\omega,w)\\
&\sim \frac{\varrho\delta_{\#}(\omega)}{2\pi\sqrt{-1}}\int_{\mathrm{Re}(s_{\tau 1})=r_1>1}\frac{\Sigma_1(R_4\Phi_{\mathfrak{p}_{41}},s_{\tau 1}+1)}{s_{\tau 1}-1}\phi_{k_1}(s_{\tau 1})\widetilde{\Lambda}(w;s_{\tau 1},1)ds_\tau.
\end{aligned}$$

*Proof.* By (6.6) and (7.5)(1), the residue of the integrand of $\Xi_{\mathfrak{p}_{41}}(\Phi,\omega,w)$ at $s_{\tau 2}=s_{\tau 1}$ is precisely $\Xi_{\mathfrak{p}_3}(\Phi,\omega,w)$. So

$$\begin{aligned}
&\Xi_{\mathfrak{p}_3}(\Phi,\omega,w)-\Xi_{\mathfrak{p}_{41}}(\Phi,\omega,w)\\
&\sim -\delta_{\#}(\omega)\left(\frac{1}{2\pi\sqrt{-1}}\right)^2\int_{\substack{\mathrm{Re}(s_\tau)=(r_1,r_2)\\ r_1>r_2>1}}\frac{\Sigma_1(R_4\Phi_{\mathfrak{p}_{41}},s_{\tau 1}+1)}{s_{\tau 2}-s_{\tau 1}}\widetilde{\Lambda}_\tau(w;s_\tau)ds_\tau.
\end{aligned}$$

If $\tau\neq(\nu,\nu)$, we can ignore the corresponding term as before. So we assume $\tau=(\nu,\nu)$. Then the above integral is

$$\begin{aligned}
&-\delta_{\#}(\omega)\left(\frac{1}{2\pi\sqrt{-1}}\right)^2\int_{\mathrm{Re}(s_\tau)=(2,\frac{1}{2})}\frac{\Sigma_1(R_4\Phi_{\mathfrak{p}_{41}},s_{\tau 1}+1)}{s_{\tau 2}-s_{\tau 1}}\widetilde{\Lambda}_\tau(w;s_\tau)ds_\tau\\
&+\frac{\varrho\delta_{\#}(\omega)}{2\pi\sqrt{-1}}\int_{\mathrm{Re}(s_{\tau 1})=r_1>1}\frac{\Sigma_1(R_4\Phi_{\mathfrak{p}_{41}},s_{\tau 1}+1)}{s_{\tau 1}-1}\phi_{k_1}(s_{\tau 1})\widetilde{\Lambda}(w;s_{\tau 1},1)ds_\tau\\
&\sim \frac{\varrho\delta_{\#}(\omega)}{2\pi\sqrt{-1}}\int_{\mathrm{Re}(s_{\tau 1})=r_1>1}\frac{\Sigma_1(R_4\Phi_{\mathfrak{p}_{41}},s_{\tau 1}+1)}{s_{\tau 1}-1}\phi_{k_1}(s_{\tau 1})\widetilde{\Lambda}(w;s_{\tau 1},1)ds_\tau.
\end{aligned}$$

Q.E.D.

Let $\tau = (\nu, \nu)$. We define

$$
\begin{aligned}
(8.2) \qquad J_1(\Phi, \omega, s_{\tau 1}) &= T^1_{\mathrm{st},1}\left(R_{0,1}\Phi, \omega_{\mathrm{st},1}, 1, -\frac{s_{\tau 1}-1}{2}\right) \\
&\quad + T^1_{\mathrm{st},2}\left(\Phi, \omega, -\frac{s_{\tau 1}-1}{2}\right) \\
&\quad + \delta_{\#}(\omega)\Sigma_1(R_4\Phi_{\mathfrak{p}_{41}}, s_{\tau 1}+1) \\
&\quad + \delta_{\#}(\omega)\Sigma_1(R_4\Phi_{\mathfrak{p}_{42}}, s_{\tau 1}+1) \\
J_2(\Phi, \omega, s_{\tau 1}) &= \frac{J_1(\Phi, \omega, s_{\tau 1})}{s_{\tau 1}-1} + \Sigma_{\mathfrak{d}_1}(\Phi, \omega, s_{\tau 1}, 1), \\
J_3(\Phi, \omega) &= \delta_{\#}(\omega)\mathfrak{V}_{k_1,2}\mathfrak{V}_2\Phi(0) + \delta(\omega_2)Z_{W,(0)}(R_2\Phi, (1, \omega_1), 2), \\
J_4(\Phi, \omega, s_\tau) &= J_2(\widehat{\Phi}, \omega^{-1}, s_{\tau 1}) - J_2(\Phi, \omega, s_{\tau 1}), \\
J_5(\Phi, \omega) &= J_3(\widehat{\Phi}, \omega^{-1}) - J_3(\Phi, \omega).
\end{aligned}
$$

Let $J_{2,(0)}(\Phi, \omega, 1)$ be the constant term of the Laurent expansion of $J_2(\Phi, \omega, s_{\tau 1})$ at $s_{\tau 1} = 1$. By (2.18), (2.22), (4.16), (5.13), (6.5), (7.5), (7.6), (8.1),

$$
\begin{aligned}
I(\Phi, \omega, w) &\sim C_G\Lambda(w; \rho)J_5(\Phi, \omega) \\
&\quad + \frac{\varrho}{2\pi\sqrt{-1}}\int_{\mathrm{Re}(s_{\tau 1})=r_1>1} J_4(\Phi, \omega, s_{\tau 1})\phi_{k_1}(s_{\tau 1})\widetilde{\Lambda}(w; s_{\tau 1}, 1)ds_{\tau 1}.
\end{aligned}
$$

Therefore, by Wright's principle (3.7.1) [5], $J_2(\Phi, \omega, s_{\tau 1})$ is holomorphic at $s_{\tau 1} = 1$ and

$$
I(\Phi, \omega, w) \sim C_G\Lambda(w; \rho)(J_5(\Phi, \omega) + J_4(\Phi, \omega, 1)).
$$

To compute $J_4(\Phi, \omega, 1)$, we just have to take the constant term of the Laurent expansion of each term of $J_4(\Phi, \omega, s_{\tau 1})$, which is $J_{2,(0)}(\widehat{\Phi}, \omega^{-1}, 1) - J_{2,(0)}(\Phi, \omega, 1)$. Let $\Sigma_{\mathfrak{d}_1,(0)}(\Phi, \omega, 1, 1)$ be the constant term of the Laurent expansion of $\Sigma_{\mathfrak{d}_1}(\Phi, \omega, s_{\tau 1}, 1)$ at $s_{\tau 1} = 1$. Then

$$
\begin{aligned}
J_{2,(0)}(\Phi, \omega, 1) = &-\frac{1}{2}T^1_{\mathrm{st},1}(\Phi, \omega, 1) - \frac{1}{2}T^1_{\mathrm{st},2}(\Phi, \omega) \\
&+ \delta_{\#}(\omega)(\Sigma_{1,(1)}(R_4\Phi_{\mathfrak{p}_{41}}, 2) + \Sigma_{1,(1)}(R_4\Phi_{\mathfrak{p}_{42}}, 2)) \\
&+ \Sigma_{\mathfrak{d}_1,(0)}(\Phi, \omega, 1, 1).
\end{aligned}
$$

It is easy to see that

$$
\begin{aligned}
\Sigma_{\mathfrak{d}_1,(0)}(\Phi, \omega, 1, 1) = &\frac{\delta_{\mathfrak{d}_1}(\omega)}{2}\Sigma_{1,1,(0,0)}(R_1\Phi, \omega_{\mathfrak{d}_1}, 1, 1) \\
&+ \frac{\delta_{\mathfrak{d}_1}(\omega)}{4}\Sigma_{1,1,(-1,1)}(R_1\Phi, \omega_{\mathfrak{d}_1}, 1, 1) \\
&+ \delta_{\mathfrak{d}_1}(\omega)\Sigma_{1,1,(1,-1)}(R_1\Phi, \omega_{\mathfrak{d}_1}, 1, 1)).
\end{aligned}
$$

By comparing the Laurent expansions of both sides of the equation

$$
\Sigma_{1,1}(R_1\Phi_\lambda, \omega_{\mathfrak{d}_1}, s_1, s_2) = \lambda^{-s_1-2s_2-1}\Sigma_{1,1}(R_1\Phi, \omega_{\mathfrak{d}_1}, s_1, s_2),
$$

we get

$$
\begin{aligned}
(8.3)\quad \Sigma_{1,1,(0,0)}(R_1\Phi_\lambda, \omega_{\mathfrak{d}_1}, 1, 1) &= \lambda^{-4}\Sigma_{1,1,(0,0)}(R_1\Phi, \omega_{\mathfrak{d}_1}, 1, 1) \\
&\quad - \lambda^{-4}\log\lambda\Sigma_{1,1,(-1,0)}(R_1\Phi, \omega_{\mathfrak{d}_1}, 1, 1) \\
&\quad - 2\lambda^{-4}\log\lambda\Sigma_{1,1,(0,-1)}(R_1\Phi, \omega_{\mathfrak{d}_1}, 1, 1) \\
&\quad + 2\lambda^{-4}(\log\lambda)^2\Sigma_{1,1,(-1,-1)}(R_1\Phi, \omega_{\mathfrak{d}_1}, 1, 1), \\
\Sigma_{1,1,(-1,1)}(R_1\Phi_\lambda, \omega_{\mathfrak{d}_1}, 1, 1) &= \lambda^{-4}\Sigma_{1,1,(-1,1)}(R_1\Phi, \omega_{\mathfrak{d}_1}, 1, 1) \\
&\quad - 2\lambda^{-4}\log\lambda\Sigma_{1,1,(-1,0)}(R_1\Phi, \omega_{\mathfrak{d}_1}, 1, 1) \\
&\quad + 2\lambda^{-4}(\log\lambda)^2\Sigma_{1,1,(-1,-1)}(R_1\Phi, \omega_{\mathfrak{d}_1}, 1, 1), \\
\Sigma_{1,1,(1,-1)}(R_1\Phi_\lambda, \omega_{\mathfrak{d}_1}, 1, 1) &= \lambda^{-4}\Sigma_{1,1,(1,-1)}(R_1\Phi, \omega_{\mathfrak{d}_1}, 1, 1) \\
&\quad - \lambda^{-4}\log\lambda\Sigma_{1,1,(0,-1)}(R_1\Phi, \omega_{\mathfrak{d}_1}, 1, 1) \\
&\quad + \frac{\lambda^{-4}(\log\lambda)^2}{2}\Sigma_{1,1,(-1,-1)}(R_1\Phi, \omega, 1, 1).
\end{aligned}
$$

We define

$$
\begin{aligned}
(8.4)\quad F_{-1}(\Phi, \omega, 4) &= \frac{\delta_{\mathfrak{d}_1}(\omega)}{2}\Sigma_{1,1,(0,0)}(R_1\Phi, \omega_{\mathfrak{d}_1}, 1, 1) \\
&\quad + \frac{\delta_{\mathfrak{d}_1}(\omega)}{4}\Sigma_{1,1,(-1,1)}(R_1\Phi, \omega_{\mathfrak{d}_1}, 1, 1) \\
&\quad + \delta_{\mathfrak{d}_1}(\omega)\Sigma_{1,1,(1,-1)}(R_1\Phi, \omega_{\mathfrak{d}_1}, 1, 1), \\
F_{-2}(\Phi, \omega, 4) &= \delta_{\mathfrak{d}_1}(\omega)\Sigma_{1,1,(-1,0)}(R_1\Phi, \omega_{\mathfrak{d}_1}, 1, 1) \\
&\quad + 2\delta_{\mathfrak{d}_1}(\omega)\Sigma_{1,1,(0,-1)}(R_1\Phi, \omega_{\mathfrak{d}_1}, 1, 1), \\
F_{-3}(\Phi, \omega, 4) &= 2\delta_{\mathfrak{d}_1}(\omega)\Sigma_{1,1,(-1,-1)}(R_1\Phi, \omega_{\mathfrak{d}_1}, 1, 1).
\end{aligned}
$$

Then

$$
\begin{aligned}
(8.5)\quad \Sigma_{\mathfrak{d}_1,(0)}(R_1\Phi_\lambda, \omega, 1, 1) &= \lambda^{-4}F_{-1}(\Phi, \omega, 4) \\
&\quad - \lambda^{-4}\log\lambda F_{-2}(\Phi, \omega, 4) \\
&\quad + \lambda^{-4}(\log\lambda)^2 F_{-3}(\Phi, \omega, 4).
\end{aligned}
$$

Note that $\Phi_{\mathfrak{p}_{41}} = R_{W,1}R_2\Phi$, $\Phi_{\mathfrak{p}_{42}} = R_{W,1}\mathscr{F}_W R_2\Phi$, and $\delta(\omega_2)\delta_{W,1}((1,\omega_1)) = \delta_\#(\omega)$.

By comparing the Laurent expansions of both sides of the equations

$$
\begin{aligned}
Z_W(R_1\Phi_\lambda, (1,\omega_1), s) &= \lambda^{-s}Z_W(R_1\Phi, (1,\omega_1), s), \\
\Sigma_1(R_{W,1}R_2\Phi_\lambda, s) &= \lambda^{-s}\Sigma_1(R_{W,1}R_2\Phi, s), \\
\Sigma_1(R_{W,1}\mathscr{F}_W R_2\Phi_\lambda, s) &= \lambda^{s-4}\Sigma_1(R_{W,1}\mathscr{F}_W R_2\Phi, s),
\end{aligned}
$$

we get

$$
\begin{aligned}
(8.6)\quad Z_{W,(0)}(R_1\Phi_\lambda, (1,\omega_1), 2) &= \lambda^{-2}Z_{W,(0)}(R_1\Phi, (1,\omega_1), 2) \\
&\quad - \lambda^{-2}\log\lambda Z_{W,(-1)}(R_1\Phi, (1,\omega_1), 2), \\
\Sigma_{1,(1)}(R_{W,1}R_2\Phi_\lambda, 2) &= \lambda^{-2}\Sigma_{1,(1)}(R_{W,1}R_2\Phi, 2) \\
&\quad - \lambda^{-2}\log\lambda\Sigma_{1,(0)}(R_{W,1}R_2\Phi, 2), \\
\Sigma_{1,(1)}(R_{W,1}\mathscr{F}_W R_2\Phi_\lambda, 2) &= \lambda^{-2}\Sigma_{1,(1)}(R_{W,1}\mathscr{F}_W R_2\Phi, 2) \\
&\quad + \lambda^{-2}\log\lambda\Sigma_{1,(0)}(R_{W,1}\mathscr{F}_W R_2\Phi, 2).
\end{aligned}
$$

We put

$$F_{-1}(\Phi,\omega,2) = \delta(\omega_2)Z_{W,(0)}(R_2\Phi,(1,\omega_1),2) + \delta_{\#}(\omega)\Sigma_{1,(1)}(R_{W,1}R_2\Phi,2) + \delta_{\#}(\omega)\Sigma_{1,(1)}(R_{W,1}\mathscr{F}_W R_2\Phi,2). \tag{8.7}$$

By (3.6),

$$\delta(\omega_2)Z_{W,(-1)}(R_2\Phi,(1,\omega_1),2) = \delta_{\#}(\omega)\Sigma_{1,(0)}(R_{W,1}\mathscr{F}_W R_2\Phi,2) + \delta_{\#}(\omega)\Sigma_{1,(0)}(R_{W,1}R_2\Phi,2).$$

So by (8.6),

$$\lambda^{-2}F_{-1}(\Phi,\omega,2) = \delta(\omega_2)Z_{W,(0)}(R_2\Phi_\lambda,(1,\omega_1),2) + \delta_{\#}(\omega)\Sigma_{1,(1)}(R_{W,1}R_2\Phi_\lambda,2) + \delta_{\#}(\omega)\Sigma_{1,(1)}(R_{W,1}\mathscr{F}_W R_2\Phi_\lambda,2). \tag{8.8}$$

It is easy to see that

$$\Phi_\lambda(0) = \Phi(0),\ \widehat{\Phi_\lambda}(0) = \lambda^{-8}\widehat{\Phi}(0). \tag{8.9}$$

Then by (8.5), (8.8), (8.9),

$$\begin{aligned} &J_3(\Phi_\lambda,\omega) + J_{2,(0)}(\Phi_\lambda,\omega,1) \\ &= -\frac{1}{2}T^1_{\mathrm{st},1}(\Phi_\lambda,\omega,1) - \frac{1}{2}T^1_{\mathrm{st},2}(\Phi_\lambda,\omega) \\ &\quad + \delta_{\#}(\omega)\mathfrak{V}_{2,k_1}\mathfrak{V}_2\Phi(0) + \lambda^{-2}F_{-1}(\Phi,\omega,2) \\ &\quad + \sum_{i=1,2,3}\lambda^{-4}(-\log\lambda)^{i-1}F_{-i}(\Phi,\omega,4). \end{aligned} \tag{8.10}$$

By similar computations, it is possible to show

$$\begin{aligned} &J_3(\widehat{\Phi_\lambda},\omega^{-1}) + J_{2,(0)}(\widehat{\Phi_\lambda},\omega^{-1},1) \\ &= -\frac{1}{2}T^1_{\mathrm{st},1}(\widehat{\Phi_\lambda},\omega^{-1},1) - \frac{1}{2}T^1_{\mathrm{st},2}(\widehat{\Phi_\lambda},\omega^{-1}) \\ &\quad + \lambda^{-8}\delta_{\#}(\omega)\mathfrak{V}_{2,k_1}\mathfrak{V}_2\widehat{\Phi}(0) + \lambda^{-6}F_{-1}(\widehat{\Phi},\omega^{-1},2) \\ &\quad + \sum_{i=1,2,3}\lambda^{-4}(\log\lambda)^{i-1}F_{-i}(\widehat{\Phi},\omega^{-1},4). \end{aligned} \tag{8.11}$$

The verification of (8.11) is left to the reader.

We put

$$F(\Phi,\omega,s) = \frac{F_{-1}(\Phi,\omega,2)}{s-2} + \sum_{i=1,2,3}\frac{F_{-i}(\Phi,\omega,4)}{(s-4)^i}. \tag{8.12}$$

**Definition (8.13)**

$$Z_{\mathrm{ad}}(\Phi,\omega,s)=Z(\Phi,\omega,s)-\frac{1}{2}T_{\mathrm{st},1}(\Phi,\omega,s,1)-\frac{1}{2}T_{\mathrm{st},2}(\Phi,\omega,s).$$

We call $Z_{\mathrm{ad}}(\Phi,\omega,s)$ the adjusted zeta function.

It is easy to see that

$$(8.14)\qquad \int_0^1 \lambda^s T^1_{\mathrm{st},1}(R_{0,1}\Phi_\lambda,\omega,1)=T_{\mathrm{st},1}(R_{0,1}\Phi,\omega,1)-T_{\mathrm{st},1+}(R_{0,1}\Phi,\omega,1),$$
$$\int_0^1 \lambda^s T^1_{\mathrm{st},1}(R_{0,1}\widehat{\Phi_\lambda},\omega,1)=T_{\mathrm{st},1+}(R_{0,1}\widehat{\Phi},\omega,1).$$

The integrals of $T^1_{\mathrm{st},2}(R_{0,2}\Phi_\lambda,\omega)$, etc. are similar.

Therefore, by (8.10, (8.11), (8.14), we get the following theorem.

**Theorem (8.15)** *Suppose that* $\Phi=M_\omega\Phi$. *Then*

$$\begin{aligned}
Z_{\mathrm{ad}}(\Phi,\omega,s)=\;&Z_+(\Phi,\omega,s)+Z_+(\widehat{\Phi},\omega^{-1},8-s)\\
&-\frac{1}{2}T_{\mathrm{st},1+}(\Phi,\omega,s,1)-\frac{1}{2}T_{\mathrm{st},1+}(\widehat{\Phi},\omega^{-1},8-s,1)\\
&-\frac{1}{2}T_{\mathrm{st},2+}(\Phi,\omega,s)-\frac{1}{2}T_{\mathrm{st},2+}(\widehat{\Phi},\omega^{-1},8-s)\\
&+\delta_\#(\omega)\mathfrak{V}_{k_1,2}\mathfrak{V}_2\left(\frac{\widehat{\Phi}(0)}{s-8}-\frac{\Phi(0)}{s}\right)+F(\widehat{\Phi},\omega^{-1},8-s)-F(\Phi,\omega,s).
\end{aligned}$$

Since $T_{\mathrm{st},1}(\Phi,\omega,s,1)$ is holomorphic for $\mathrm{Re}(s)>5$ and $T_{\mathrm{st},2}(\Phi,\omega,s)$ is holomorphic for $\mathrm{Re}(s)>6$, we get the following corollary.

**Corollary (8.16)** *The unadjusted zeta function* $Z(\Phi,\omega,s)$ *can be continued meromorphically to* $\mathrm{Re}(s)>6$ *having the only possible simple pole at* $s=8$ *with residue* $\delta_\#(\omega)\mathfrak{V}_{k_1,2}\mathfrak{V}_2\widehat{\Phi}(0)$.

**Remark (8.17)** Note that

$$F_{-3}(\Phi,\omega,4)=\delta_\#(\omega)\int_{\mathbb{A}_{k_1}\times(\mathbb{A})^2}\widetilde{R}_1\Phi(x_{12},x_{21},x_{22})dx_{21}dx_{12}dx_{22}.$$

If $\Phi=M_\omega\Phi$ and $\delta_\#(\omega)=1$, $\Phi$ is $K$-invariant. Then it is easy to see that

$$F_{-3}(\Phi,\omega,4)=F_{-3}(\widehat{\Phi},\omega^{-1},4).$$

Therefore, the pole at $s=4$ is a double pole. So the orders of the poles of $Z_{\mathrm{ad}}(\Phi,\omega,s)$ at $s=0,2,4,6,8$ are $1,1,2,1,1$ respectively.

Akihiko Yukie
Oklahoma State University
Mathematics Department
401 Math Science
Stillwater OK 74078-1058 USA
yukie@math.okstate.edu